\documentclass[11pt,a4paper,reqno]{amsart}

\usepackage[utf8]{inputenc}
\usepackage[T1]{fontenc}
\usepackage{lmodern}
\usepackage{amssymb}
\usepackage[a4paper,margin=1.2in]{geometry}
\usepackage[colorlinks=true,linkcolor=blue,citecolor=blue,urlcolor=blue]{hyperref}

\theoremstyle{plain}
\newtheorem{theorem}{Theorem}[section]
\newtheorem{lemma}[theorem]{Lemma}
\newtheorem{proposition}[theorem]{Proposition}
\newtheorem{corollary}[theorem]{Corollary}

\theoremstyle{definition}
\newtheorem{definition}[theorem]{Definition}

\theoremstyle{remark}
\newtheorem{remark}[theorem]{Remark}

\numberwithin{equation}{section}

\newcommand{\RR}{\mathbb{R}}
\newcommand{\CC}{\mathbb{C}}
\newcommand{\NN}{\mathbb{N}}
\newcommand{\fra}{\mathfrak{a}}
\newcommand{\frb}{\mathfrak{b}}
\newcommand{\frc}{\mathfrak{c}}
\newcommand{\frn}{\mathfrak{n}}
\newcommand{\cX}{\mathcal{X}}
\newcommand{\cF}{\mathcal{F}}
\newcommand{\cK}{\mathcal{K}}
\newcommand{\cP}{\mathcal{P}}
\newcommand{\cS}{\mathcal{S}}
\newcommand{\dd}{\mathrm{d}}
\newcommand{\tq}{\widetilde{q}}
\newcommand{\tA}{\widetilde{A}}
\newcommand{\tla}{\widetilde{\lambda}}
\newcommand{\tphi}{\widetilde{\phi}}
\newcommand{\tpsi}{\widetilde{\psi}}
\newcommand{\tPsi}{\widetilde{\Psi}}
\newcommand{\ttheta}{\widetilde{\theta}}
\newcommand{\tm}{\widetilde{m}}
\newcommand{\ta}{\widetilde{a}}
\newcommand{\tc}{\widetilde{c}}
\newcommand{\tG}{\widetilde{G}}
\newcommand{\tM}{\widetilde{M}}
\newcommand{\tH}{\widetilde{H}}
\newcommand{\tK}{\widetilde{K}}
\newcommand{\abs}[1]{\lvert #1\rvert}
\newcommand{\norm}[1]{\lVert #1\rVert}
\newcommand{\ps}[2]{\langle #1,#2\rangle}

\DeclareMathOperator{\supp}{supp}

\DeclareMathOperator{\re}{Re}

\begin{document}
	
	\title[Fixed-time inverse problem for the Robin Laplacian]
	{Determination of the unbounded potential of the Robin Laplacian
		from fixed-time boundary measurements of time-fractional
		diffusion-wave equations}
	
	\author{Abdelmalek Metidji}
	\address{Aix Marseille University, CNRS, CPT,
		Marseille, France}
	\email{abdelmalek.metidji1@gmail.com}
	\email{abdelmalek.metidji@univ-amu.fr}
	
	\author{Yacine Mokhtari}
	\address{Research Center in Applied Mathematics, Sidi Abdellah Scientific Hub, Mahelma, Algiers, Algeria}
	\email{mokhtari.yacine@crma.dz}
	
	
	\begin{abstract}
		We consider the time-fractional diffusion-wave equation
		$\partial_t^\rho u - \Delta u + q u = 0$ in a bounded domain
		$\Omega \subset \RR^n$, $n \geq 3$, where $\partial_t^\rho$ is the Caputo
		derivative of arbitrary order $\rho \in (0,2)$. This framework covers in a
		unified way the sub-diffusive case $\rho \in (0,1)$, the classical heat
		equation $\rho = 1$, and the super-diffusive case $\rho \in (1,2)$.
		The real-valued potential $q$ is only assumed to lie in
		$L^{n/2}(\Omega)$, and the equation is endowed with zero initial
		conditions and a non-homogeneous Robin boundary condition, whose
		coefficient is allowed to be unbounded as well. The boundary excitations
		are supported in an arbitrary open subset $\Gamma_0$ of the boundary,
		and the trace of the solution is measured, at one fixed time, on another
		open subset $\Gamma_1$. We prove that the corresponding
		Robin-to-Dirichlet map uniquely determines the boundary spectral data
		of the underlying Robin Laplacian on $\Gamma_0 \cup \Gamma_1$,
		provided that $\Gamma_0$ and $\Gamma_1$ overlap. When
		$\Gamma_0 \cup \Gamma_1$ is the whole boundary, we deduce that $q$
		is uniquely determined by this map in every dimension $n \geq 3$,
		the critical case $q \in L^{3/2}(\Omega)$ in dimension three being
		included. The proof combines analytic continuation in the time variable,
		the Laplace transform and an identification of the poles and residues
		of meromorphic functions. The potential is then recovered from the
		boundary spectral data through an elliptic Robin-to-Dirichlet map,
		an integral identity, and complex geometric optics solutions.
	\end{abstract}
	
	\maketitle
\section{Introduction}\label{sec:intro}

\subsection{Statement of the problem}\label{sec:problem}

Let $\Omega$ be a bounded domain of $\RR^n$, $n \geq 3$, with $C^{1,1}$ boundary
$\Gamma$, and let $\nu$ denote the outward unit normal vector field to $\Gamma$.
Given $T \in (0,\infty)$, we set $Q := \Omega \times (0,T)$ and
$\Sigma := \Gamma \times (0,T)$. Let $\rho \in (0,2)$ and let $m \in \{1,2\}$ be
the smallest integer such that $\rho \leq m$. We consider the initial-boundary
value problem (IBVP)
\begin{equation}\label{eq:ibvp}
\begin{cases}
\partial_t^\rho u - \Delta u + q u = 0 & \text{in } Q,\\
(\partial_\nu + \alpha) u = g & \text{on } \Sigma,\\
\partial_t^j u(\cdot,0) = 0 \ \text{ in } \Omega, & j = 0,\ldots,m-1,
\end{cases}
\end{equation}
where the real-valued functions $q : \Omega \to \RR$ and
$\alpha : \Gamma \to \RR$ do not depend on the time variable $t$, and where
$\partial_t^\rho$ denotes the Caputo fractional derivative of order $\rho$ with
respect to $t$. Namely, $\partial_t^1 := \partial_t$ and, for
$\rho \in (0,1) \cup (1,2)$,
\begin{equation}\label{eq:caputo}
\partial_t^\rho y(t) := \frac{1}{\Gamma(m-\rho)} \int_0^t (t-\tau)^{m-\rho-1}
y^{(m)}(\tau)\, \dd\tau, \quad t \in (0,T],
\end{equation}
for every function $y \in C^m([0,T])$. Here and in the remaining part of this
text, $\Gamma(\cdot)$ stands for the Euler Gamma function; since it always
appears with an argument, it cannot be confused with the boundary $\Gamma$ of
$\Omega$.

When $\rho = 1$, the first line of \eqref{eq:ibvp} is the classical heat
equation, which models, for instance, heat conduction in $\Omega$. In the
sub-diffusive case $\rho \in (0,1)$, it is usually called the fractional
diffusion equation, while in the case $\rho \in (1,2)$ it is often referred to
as the fractional wave equation. These fractional equations describe anomalous
diffusion phenomena, such as diffusion in heterogeneous or porous media, see
e.g.\ \cite{MK00}. The Robin boundary
condition models the exchange of the diffusing quantity with the surroundings of
$\Omega$, $\alpha$ being the transfer coefficient and $g$ an imposed boundary
flux, while $q$ accounts for absorption. In other words, $u$ is the response of
the system \eqref{eq:ibvp} to the boundary excitation $g$.

In this article we are concerned with the inverse problem of determining the
potential $q$ from boundary measurements of $u$ performed at one fixed time, the
order $\rho$ and the Robin coefficient $\alpha$ being known. More precisely, let
$\Gamma_0$ and $\Gamma_1$ be two nonempty open subsets of $\Gamma$ and let
$T_0 \in (0,T]$. We probe the system with excitations $g$ supported in
$\Gamma_0 \times (0,T_0)$ and we observe the trace on $\Gamma_1$ of the solution
$u$ to \eqref{eq:ibvp} at time $T_0$. The data of the inverse problem are thus
encoded in the operator
\[
D_q : g \mapsto u(\cdot,T_0)_{|\Gamma_1},
\]
which we call the \emph{Robin-to-Dirichlet map} at time $T_0$; its
rigorous definition is given in Section~\ref{sec:RtD} below. Our goal is to show
that the mapping $q \mapsto D_q$ is injective. We stress that neither the
potential $q$ nor the Robin coefficient $\alpha$ is assumed to be bounded, that
no sign condition is imposed on $q$, and that $\alpha$ may change sign, as it is
only required to be bounded from below by a suitable negative constant, see
\eqref{eq:alpha} below.

\subsection{The Robin Laplacian}\label{sec:robin}

We stick to the framework of \cite{CMS24}, which we briefly recall for the
reader's convenience. Let $H := L^2(\Omega)$ and $V := H^1(\Omega)$ be endowed
with their usual scalar products and norms, the scalar product of $H$ being
denoted by $(\cdot,\cdot)$. Put $p := 2n/(n+2)$ and let $p^* := 2n/(n-2)$ be its
conjugate exponent, in such a way that $V$ is continuously embedded in
$L^{p^*}(\Omega)$. For any measurable subset $\Gamma'$ of $\Gamma$, we write
\[
\ps{f}{h}_{\Gamma'} := \int_{\Gamma'} f\, \overline{h}\, \dd\sigma,
\quad f,\, h \in L^2(\Gamma'),
\]
where $\dd\sigma$ denotes the surface measure on $\Gamma$, and we denote by
$\frn$ the norm of the (bounded) trace operator
$u \in V \mapsto u_{|\Gamma} \in L^2(\Gamma)$, so that
\begin{equation}\label{eq:trace}
\norm{u_{|\Gamma}}_{L^2(\Gamma)} \leq \frn \norm{u}_V, \quad u \in V.
\end{equation}

Let $s \in (n-1,\infty)$. For $\alpha \in L^s(\Gamma,\RR)$ and
$q \in L^{n/2}(\Omega,\RR)$, we consider the sesquilinear form
\begin{equation}\label{eq:form}
\fra(u,v) := \int_\Omega \nabla u \cdot \nabla \overline{v}\, \dd x
+ \int_\Omega q u \overline{v}\, \dd x
+ \int_\Gamma \alpha u \overline{v}\, \dd\sigma, \quad u,\, v \in V,
\end{equation}
which is continuous on $V \times V$, see \cite[Section~1.1]{CMS24}. We denote by
$C_\fra$ a positive constant such that
\begin{equation}\label{eq:acont}
\abs{\fra(u,v)} \leq C_\fra \norm{u}_V \norm{v}_V, \quad u,\, v \in V.
\end{equation}
Throughout the entire text we assume, as in \cite{CMS24}, that
\begin{equation}\label{eq:alpha}
\alpha \geq -\frc\ \text{ a.e.\ on } \Gamma,
\ \text{ for some constant } \frc \in (0,\frn^{-2}).
\end{equation}
This condition is obviously satisfied when $\alpha$ is non-negative. With
reference to \cite[Eq.~(1.2)]{CMS24}, there exist two constants
$\kappa \in (0,1)$ and $\lambda^* \in (0,\infty)$, depending on $\Omega$,
$\alpha$ and $q$, such that
\begin{equation}\label{eq:coer}
\fra(u,u) + \lambda^* \norm{u}_H^2 \geq \kappa \norm{u}_V^2, \quad u \in V.
\end{equation}

Since $q$ and $\alpha$ are real-valued, the form $\fra$ is Hermitian. Hence the
bounded operator $A : V \to V^*$, where $V^*$ is the dual space of $V$, defined
by $\ps{Au}{v}_{V^*,V} = \fra(u,v)$ for all $u,\, v \in V$, is self-adjoint.
Notice that $A$ itself need not be coercive, since $\fra(u,u)$ may be negative:
for instance, if $\alpha = 0$ and $q = -1$, then $\fra(1,1) = -\abs{\Omega}$.
Nevertheless, upon identifying $H$ with a subspace of $V^*$ in the usual way, the
operator $A + \lambda^*$ is coercive by \eqref{eq:coer}. Consequently, since $V$
is compactly embedded in $H$, we infer from \cite[Theorem~2.37]{McL00}, applied
to $A + \lambda^*$ (see also \cite[Section~1.2]{CMS24}), that the spectrum of $A$
consists of its eigenvalues
$\lambda_k$, $k \in \NN := \{1,2,\ldots\}$, arranged in non-decreasing order and
repeated with the (finite) multiplicity,
\[
-\lambda^* < \lambda_1 \leq \lambda_2 \leq \cdots \leq \lambda_k \leq \cdots,
\qquad \lim_{k \to \infty} \lambda_k = +\infty,
\]
and there exists an orthonormal basis $\{\phi_k,\ k \in \NN\}$ of $H$, made of
eigenfunctions $\phi_k \in V$ of $A$, i.e.
\begin{equation}\label{eq:eig}
\fra(\phi_k,v) = \lambda_k (\phi_k,v), \quad v \in V,\ k \in \NN.
\end{equation}
Since $q$ and $\alpha$ are real-valued, $\overline{\phi}$ solves
\eqref{eq:eig} whenever $\phi$ does. Therefore each eigenspace of $A$ is spanned
by real-valued functions, and the $\phi_k$ may be, and always will be, chosen
real-valued. For the sake of shortness we write
\[
\psi_k := {\phi_k}_{|\Gamma}, \quad k \in \NN,
\]
and we recall that the boundary spectral data (BSD) of $A$ are defined as
\[
\mathrm{BSD}(q) := \{ (\lambda_k,\psi_k),\ k \in \NN \}.
\]
Notice that the $\psi_k$, and hence the BSD, depend on the choice of the
orthonormal basis $\{\phi_k,\ k \in \NN\}$.

\subsection{Weak solutions and the Robin-to-Dirichlet map}\label{sec:RtD}

For $\beta \in (0,\infty)$ we put
\begin{equation}\label{eq:kbeta}
k_\beta(t) := \frac{t^{\beta-1}}{\Gamma(\beta)}, \quad t \in (0,\infty),
\end{equation}
and we denote by $I^\beta$ the Riemann--Liouville fractional integral of order
$\beta$, i.e.
\begin{equation}\label{eq:RL}
I^\beta w(t) := (k_\beta * w)(t) := \int_0^t k_\beta(t-\tau)\, w(\tau)\,
\dd\tau, \quad w \in L^1(0,T),\ \text{a.e. } t \in (0,T).
\end{equation}

Assume for a while that $u$ is a sufficiently smooth solution to
\eqref{eq:ibvp}. Multiplying the first line of \eqref{eq:ibvp} by
$\overline{v}$, where $v \in V$, integrating by parts over $\Omega$ and taking
the second line of \eqref{eq:ibvp} into account, we get that
\[
\partial_t^\rho (u(\cdot,t),v) + \fra(u(\cdot,t),v)
= \ps{g(\cdot,t)}{v_{|\Gamma}}_\Gamma, \quad t \in (0,T).
\]
Applying $I^\rho$ to both sides of this identity and taking into account the
initial conditions of \eqref{eq:ibvp}, which ensure that
$I^\rho \partial_t^\rho y = y$ for $y(t) := (u(\cdot,t),v)$ (see
Remark~\ref{rem:caputo} below), we end up getting the integral identity
\eqref{eq:weak} below. This leads us to the following definition, where
$L^2(\Sigma)$ is identified with $L^2(0,T;L^2(\Gamma))$.

\begin{definition}\label{def:weak}
Let $g \in L^2(\Sigma)$. A \emph{weak solution} to the IBVP \eqref{eq:ibvp} is a
function $u \in L^2(0,T;V)$ such that, for every $v \in V$, we have
\begin{equation}\label{eq:weak}
(u(\cdot,t),v) + I^\rho\big[ \fra(u,v) \big](t)
= I^\rho\big[ \ps{g}{v_{|\Gamma}}_\Gamma \big](t)
\quad \text{for a.e. } t \in (0,T),
\end{equation}
where $\fra(u,v)$ denotes the function $t \mapsto \fra(u(\cdot,t),v)$ and
$\ps{g}{v_{|\Gamma}}_\Gamma$ denotes the function
$t \mapsto \ps{g(\cdot,t)}{v_{|\Gamma}}_\Gamma$; both of them belong to
$L^2(0,T)$.
\end{definition}

When $\rho = 1$, \eqref{eq:weak} is the time-integrated form of the usual
variational formulation of the IBVP associated with the heat equation. We shall
see in
Proposition~\ref{prop:wp} that for every $g \in L^2(\Sigma)$, the IBVP
\eqref{eq:ibvp} admits a unique weak solution, which we denote by $u_g$.

Notice that the initial conditions of \eqref{eq:ibvp} are encoded in the
integral identity \eqref{eq:weak}, and that they are understood in this weak
sense only. Indeed, for a general excitation $g \in L^2(\Sigma)$, the weak
solution $u_g$ need not have a trace at $t = 0$. For instance, let $q = 0$,
$\alpha = 0$, $\rho \in (0,1/2)$, $\beta \in (\rho,1/2)$ and
$g(\sigma,t) := t^{-\beta}$, which belongs to $L^2(\Sigma)$ since $\beta < 1/2$.
Then, taking $v = 1$ in \eqref{eq:weak}, we get that
\[
(u_g(\cdot,t),1) = \abs{\Gamma}\, I^\rho\big[ \tau \mapsto \tau^{-\beta} \big](t)
= \abs{\Gamma}\, \frac{\Gamma(1-\beta)}{\Gamma(1+\rho-\beta)}\, t^{\rho-\beta}
\quad \text{for a.e. } t \in (0,T),
\]
where $\abs{\Gamma}$ denotes the surface measure of $\Gamma$, and the right-hand
side is unbounded near $t = 0$, since $\rho < \beta$. In contrast, if $g$
vanishes a.e.\ in $\Gamma \times (0,\varepsilon)$ for some $\varepsilon > 0$,
which is the case of the excitations of the set $\cX_0$ defined below, then
$u_g$ vanishes a.e.\ in $\Omega \times (0,\varepsilon)$, as can be seen from the
representation formula \eqref{eq:rep} below.

Next, we define the set of admissible excitations as
\begin{equation}\label{eq:X0}
\cX_0 := \{ g \in L^2(\Sigma);\ \supp g \subset \Gamma_0 \times (0,T_0) \},
\end{equation}
where $\supp g$ is the support of $g$, i.e.\ the complement in
$\Gamma \times [0,T]$ of the largest relatively open subset on which $g$
vanishes almost everywhere. As $\supp g$ is compact, every $g \in \cX_0$
vanishes a.e.\ in $\Gamma \times \big( (0,\varepsilon) \cup
(T_0-\varepsilon,T) \big)$ for some
$\varepsilon = \varepsilon(g) \in (0,T_0/2)$. We shall see in
Proposition~\ref{prop:T0} that this property guarantees that $u_g$ coincides
a.e.\ on $[T_0 - \varepsilon/2,T_0]$ with a continuous function from
$[T_0 - \varepsilon/2,T_0]$ into $V$. Denoting by $u_g(\cdot,T_0) \in V$ the
value of this function at time $T_0$, the Robin-to-Dirichlet map at
time $T_0$,
\begin{equation}\label{eq:Dq}
D_q : \cX_0 \to L^2(\Gamma_1), \qquad g \mapsto u_g(\cdot,T_0)_{|\Gamma_1},
\end{equation}
is thus well defined. Notice that $D_q$ depends on $\rho$, $\alpha$, $T_0$,
$\Gamma_0$ and $\Gamma_1$, which are fixed throughout the text.

In the sequel, we write $\tA$ (resp., $\tla_k$, $\tphi_k$, $\tpsi_k$, $D_{\tq}$)
instead of $A$ (resp., $\lambda_k$, $\phi_k$, $\psi_k$, $D_q$) when the
potential $\tq$ is substituted for $q$.

\subsection{Main results}\label{sec:results}

The first result of this article states that the potential is uniquely
determined by the Robin-to-Dirichlet map, provided the excitation
zone $\Gamma_0$ and the observation zone $\Gamma_1$ overlap and cover the whole
boundary.

\begin{theorem}\label{thm:uniq}
Let $n \geq 3$, $\rho \in (0,2)$ and $T_0 \in (0,T]$. Let $\Gamma_0$ and
$\Gamma_1$ be two nonempty open subsets of $\Gamma$ such that
\[
\Gamma_0 \cup \Gamma_1 = \Gamma \quad \text{and} \quad
\Gamma_0 \cap \Gamma_1 \neq \emptyset,
\]
and let $q$ and $\tq$ be in $L^{n/2}(\Omega,\RR)$. Then we have the implication
\[
\big( D_q(g) = D_{\tq}(g) \ \text{ for all } g \in \cX_0 \big)
\Longrightarrow \big( q = \tq \ \text{ in } \Omega \big).
\]
\end{theorem}

The proof of Theorem~\ref{thm:uniq} relies on the following result, which
establishes that the Robin-to-Dirichlet map determines the boundary spectral
data on $\Gamma_0 \cup \Gamma_1$, and which does not require that
$\Gamma_0 \cup \Gamma_1$ be the whole boundary.

\begin{theorem}\label{thm:bsd}
Let $n \geq 3$, $\rho \in (0,2)$ and $T_0 \in (0,T]$. Let $\Gamma_0$ and
$\Gamma_1$ be two nonempty open subsets of $\Gamma$ such that
$\Gamma_0 \cap \Gamma_1 \neq \emptyset$, and let $q$ and $\tq$ be in
$L^{n/2}(\Omega,\RR)$. Assume that
\[
D_q(g) = D_{\tq}(g), \quad g \in \cX_0.
\]
Then we have $\lambda_k = \tla_k$ for all $k \in \NN$, and there exists an
orthonormal basis $\{ \tphi_k,\ k \in \NN \}$ of $H$, made of real-valued
eigenfunctions of $\tA$ associated with the eigenvalues $\tla_k$, such that
\[
\psi_k = \tpsi_k \ \text{ a.e.\ on } \Gamma_0 \cup \Gamma_1, \quad k \in \NN.
\]
\end{theorem}

We now collect several comments on Theorems~\ref{thm:uniq} and \ref{thm:bsd}.

\begin{remark}\label{rem:main}
\begin{enumerate}
\item[\textup{(i)}] When $\Gamma$ is connected, the condition
$\Gamma_0 \cap \Gamma_1 \neq \emptyset$ is automatically fulfilled by any two
nonempty open subsets $\Gamma_0$ and $\Gamma_1$ of $\Gamma$ such that
$\Gamma_0 \cup \Gamma_1 = \Gamma$.
\item[\textup{(ii)}] The proof of Theorem~\ref{thm:bsd} only involves
excitations with separated variables. More precisely, let $\cF$ be a set of
functions $f \in L^2(\Gamma)$ with $\supp f \subset \Gamma_0$, which is dense in
$L^2(\Gamma_0)$. Then the statements of Theorems~\ref{thm:uniq} and
\ref{thm:bsd} remain valid when the identity $D_q(g) = D_{\tq}(g)$ is only
assumed for all the excitations of the form $g(\sigma,t) = f(\sigma) h(t)$, with
$f \in \cF$ and $h \in C_c^\infty(I)$, where $I$ is an arbitrary nonempty open
interval contained in $(0,T_0)$, fixed once and for all. In other words, the
excitations may be switched on during an arbitrarily short time window only.
\item[\textup{(iii)}] No assumption other than \eqref{eq:alpha} is imposed on
the Robin coefficient $\alpha \in L^s(\Gamma,\RR)$, and no sign condition is
imposed on the potentials. In particular, the operator $A$ may have finitely
many non-positive eigenvalues. Moreover, Theorem~\ref{thm:uniq} holds for all
potentials in $L^{n/2}(\Omega,\RR)$, including the critical case
$q \in L^{3/2}(\Omega,\RR)$ when $n = 3$. When $\Gamma_0 \cup \Gamma_1 = \Gamma$,
Theorem~\ref{thm:bsd} provides the full boundary spectral data, from which $q$
could be retrieved by means of the Borg--Levinson theorem
\cite[Theorem~1.1]{CMS24} for the Robin Laplacian; but the latter is stated for
potentials in $L^r(\Omega,\RR)$ with $r > 3/2$ when $n = 3$. Instead, following
the strategy of \cite{NSU88} (see also \cite[Theorem~1.1]{Poh18} for Dirichlet
boundary conditions), we show that the boundary spectral data determine an
elliptic Robin-to-Dirichlet map, we derive from it an integral identity, and we
conclude with the complex geometric optics solutions for potentials in
$L^{n/2}$ constructed in \cite[Section~4]{KU16}, see
Section~\ref{sec:proof-uniq}.
\item[\textup{(iv)}] When $\Gamma_0 \cup \Gamma_1$ is a strict subset of
$\Gamma$, Theorem~\ref{thm:bsd} shows that $D_q$ still determines the eigenvalues
of $A$ and the traces of its eigenfunctions on $\Gamma_0 \cup \Gamma_1$. As a
consequence, $q$ is uniquely determined by $D_q$ whenever it is uniquely
determined by these partial boundary spectral data. We refer the reader to
\cite[Corollary~1.5]{CMS25} for a result in this direction, which holds for
potentials that are a priori known in a neighborhood of $\Gamma$, under
additional regularity, geometric and sign assumptions on $\Omega$, $\alpha$ and
the potentials.
\item[\textup{(v)}] Theorems~\ref{thm:uniq} and \ref{thm:bsd} are stated for all
$\rho \in (0,2)$, and they cover in a unified way the sub-diffusive case
$\rho \in (0,1)$, the classical heat equation $\rho = 1$, and the
super-diffusive case $\rho \in (1,2)$, with the same proof. When $\rho = 1$, the
Mittag-Leffler function $E_{1,1}$ used below is the exponential function, so
that $t^{\rho-1} E_{\rho,\rho}(-\lambda_k t^\rho) = e^{-\lambda_k t}$. In this
case, Theorems~\ref{thm:uniq} and \ref{thm:bsd} are the Robin counterparts of
the results of \cite{CK01,Ch24b}, obtained for Dirichlet boundary conditions.
\end{enumerate}
\end{remark}

\subsection{Comments and existing literature}\label{sec:biblio}

There is a wide mathematical literature dealing with inverse coefficient
problems for parabolic equations; we refer the reader to the monographs
\cite{Ch09,Is06} for an overview. The idea of recovering the boundary spectral
data of the underlying elliptic operator from boundary measurements of the
solution of a heat equation, performed at one fixed time, goes back to Canuto
and Kavian. In \cite{CK01} they proved that the parabolic Dirichlet-to-Neumann
map at a fixed time, where the Dirichlet data are imposed on one part of the
boundary and the flux is measured on another part, uniquely determines the
boundary spectral data of the associated Dirichlet operator, and they derived
several identification results from this; see also \cite{KKLM04} for the
equivalence between various time-domain inverse problems and inverse boundary
spectral problems. More recently, Choulli proved in \cite{Ch24b} that the
zero-order coefficient of the heat equation with Dirichlet boundary condition is
uniquely determined in the class $L^{n/2}(\Omega)$ by the parabolic
Dirichlet-to-Neumann map at a fixed time; this result is stated for $n \geq 5$
and extends to lower dimensions up to suitable modifications.

The well-posedness of initial-boundary value problems for time-fractional
diffusion-wave equations was examined in \cite{SY11}, and we refer the reader to
\cite{KRY20} for a theoretical introduction to these equations, and to
\cite{GKMR20,KST06,Pod99} for fractional calculus and Mittag-Leffler functions.
Inverse coefficient problems for such equations, from boundary measurements,
were considered in several works. In \cite{LIY16}, the number and the orders
of the fractional time derivatives of a multi-term equation, together with
space-dependent coefficients, are uniquely determined by a
Dirichlet-to-Neumann map associated with Dirichlet excitations of the form
$\lambda(t) f(x)$, where $\lambda$ is a fixed function of the time variable.
In \cite{KOSY18}, the approach of \cite{CK01} was adapted
to time-fractional diffusion-wave equations of order $\rho \in (0,1) \cup (1,2)$
with Dirichlet boundary condition: it is proved that the partial
Dirichlet-to-Neumann map at a fixed time determines two of the three bounded
coefficients (density, conductivity and non-negative potential) appearing in the
equation, as well as a smooth Riemannian manifold and smooth coefficients, up
to isometry and to the natural gauge invariance. In \cite{KLLY21}, at most two
out of four space-dependent coefficients, possibly together with the fractional
order, are uniquely recovered from a single measurement of the partial Neumann
output associated with a suitably chosen Dirichlet input; the proof relies on
the time-analyticity of the solution and on Dirichlet-to-Neumann maps in the
frequency domain.

The results of this article are also closely related to the multidimensional
Borg--Levinson theorem, which states that the boundary spectral data uniquely
determine the potential. The one-dimensional theory goes back to
\cite{Borg46,Lev49,GL51}. In higher dimension, the unique determination of a
bounded potential by the full boundary spectral data was proved by Nachman,
Sylvester and Uhlmann \cite{NSU88} and, independently, by Novikov \cite{Nov88}.
Isozaki \cite{Iso91} showed that finitely many spectral data may be omitted, and
the stability issue was addressed in \cite{AS90,CS13,KKS15}; see also the survey
\cite{Soc20}. Unbounded potentials were considered in \cite{PS02,Poh18}, and
then in \cite{BCKPS,BKMS24,Ch24a,KS22}. In the full data case, the
Borg--Levinson problem can be reduced to the determination of the potential by
an elliptic Dirichlet-to-Neumann map, as in \cite{NSU88,Poh18}; for potentials
in $L^{n/2}(\Omega)$, $n \geq 3$, the latter problem was solved by Lavine and
Nachman (see the announcement in \cite{Nac92}), and we refer the reader to
\cite{DKS13} for admissible Riemannian manifolds and to \cite{KU16} for
polyharmonic operators. All these results were obtained for
Dirichlet boundary conditions, except for \cite{NSU88}, where bounded potentials
are identified from the full boundary spectral data of the Robin Laplacian. The
case of incomplete data for the Robin Laplacian was heuristically discussed in
\cite{Smi15} for bounded potentials. It was rigorously treated in \cite{CMS24},
together with the stability issue, for unbounded potentials, and the case of
asymptotic or incomplete data was examined in \cite{CMS25}.

To the best of our knowledge, the determination of an unbounded potential by
boundary measurements performed at a fixed time has not been addressed so far
for Robin boundary conditions, neither for the heat equation nor for
time-fractional diffusion-wave equations. The present article fills this gap
in a unified way for all orders $\rho \in (0,2)$, and it can be seen as the
counterpart, for the Robin Laplacian, of the uniqueness results of
\cite{CK01,Ch24b} when $\rho = 1$, and of \cite{KOSY18} when $\rho \neq 1$.
Although the general strategy follows the one of \cite{CK01,KOSY18}, the present
framework requires several specific tools. First, since the non-homogeneous
Robin condition is a natural boundary condition, the direct problem is solved in
the weak sense of Definition~\ref{def:weak} for all excitations in
$L^2(\Sigma)$, without any lifting of the boundary data, and under the sole
assumptions that $q \in L^{n/2}(\Omega,\RR)$ and that $\alpha \in L^s(\Gamma,\RR)$
satisfies \eqref{eq:alpha}.
The weak formulation \eqref{eq:weak}, based on the Riemann--Liouville integral
$I^\rho$, avoids any regularity requirement in time. Second, the eigenfunctions
of $A$ are only known to belong to $H^1(\Omega)$, so the unique continuation
argument yielding the linear independence of their boundary traces is carried
out at the level of the weak formulation, see Lemma~\ref{lem:ucp}. Third, when
$\rho \neq 1$, the time-dependent factors $t^{\rho-1} E_{\rho,\rho}(-\lambda_k
t^\rho)$ appearing in the representation of the solution decay only like
$\lambda_k^{-2}$ as $k \to \infty$, see Lemma~\ref{lem:ML}, and the estimates
at our disposal on the traces $\psi_k$ do not guarantee that the boundary kernel
$\sum_{k \geq 1} t^{\rho-1} E_{\rho,\rho}(-\lambda_k t^\rho)\, \psi_k \otimes
\psi_k$ converges in $L^2(\Gamma \times \Gamma)$. For this reason we work with
scalar quantities, see \eqref{eq:G}, whose convergence follows from a Bessel
type inequality (Lemma~\ref{lem:spec}). Finally, the boundary spectral data are
here made of Dirichlet traces of the eigenfunctions. We show that, when
$\Gamma_0 \cup \Gamma_1 = \Gamma$, they determine the elliptic
Robin-to-Dirichlet map at a fixed energy, which yields an integral identity for
the solutions of the elliptic equations associated with the two potentials. The
potential is then recovered with the aid of the complex geometric optics
solutions of \cite{KU16}. In contrast with the Dirichlet case
\cite[Lemma~3.4]{Poh18}, the series representing the elliptic
Robin-to-Dirichlet map converges directly, by the Bessel type inequality of
Lemma~\ref{lem:spec}, and this approach allows us to deal with all potentials
in $L^{n/2}(\Omega,\RR)$, including the critical case $n = 3$.

\subsection{Outline}\label{sec:outline}

In Section~\ref{sec:prelim} we gather several properties of the eigenpairs of
the Robin Laplacian, of the Riemann--Liouville fractional integrals and of the
Mittag-Leffler functions, that are needed in the sequel. Section~\ref{sec:direct}
is devoted to the analysis of the direct problem: we prove that the IBVP
\eqref{eq:ibvp} is well posed in the sense of Definition~\ref{def:weak}, and we
establish a representation formula for the Robin-to-Dirichlet map.
Finally, Section~\ref{sec:proof} contains the proof of Theorems~\ref{thm:uniq}
and \ref{thm:bsd}.

\section{Preliminaries}\label{sec:prelim}

In Sections~\ref{sec:prelim} and \ref{sec:direct}, $q \in L^{n/2}(\Omega,\RR)$
and $\rho \in (0,2)$ are fixed.

\subsection{Spectral preliminaries}

We set
\[
\mu_k := \lambda_k + \lambda^* > 0, \quad k \in \NN,
\]
and we introduce the sesquilinear form
\[
\frb(u,v) := \fra(u,v) + \lambda^* (u,v), \quad u,\, v \in V.
\]
Since $\fra$ is Hermitian, we infer from \eqref{eq:acont} and \eqref{eq:coer}
that $\frb$ is a scalar product on $V$, which satisfies
\begin{equation}\label{eq:b}
\kappa \norm{u}_V^2 \leq \frb(u,u) \leq (C_\fra + \lambda^*) \norm{u}_V^2,
\quad u \in V.
\end{equation}
In other words, $V$ endowed with the scalar product $\frb$ is a Hilbert space,
whose norm is equivalent to the usual one. Moreover, \eqref{eq:eig} can be
equivalently rewritten as
\begin{equation}\label{eq:beig}
\frb(\phi_k,v) = \mu_k (\phi_k,v), \quad v \in V,\ k \in \NN.
\end{equation}

\begin{lemma}\label{lem:spec}
The following statements hold true.
\begin{enumerate}
\item[\textup{(i)}] A function $u \in H$ belongs to $V$ if and only if
$\sum_{k \geq 1} \mu_k \abs{(u,\phi_k)}^2 < \infty$. In this case, the series
$\sum_{k \geq 1} (u,\phi_k) \phi_k$ converges to $u$ in $V$, and we have
\begin{equation}\label{eq:parseval}
\frb(u,u) = \sum_{k \geq 1} \mu_k \abs{(u,\phi_k)}^2.
\end{equation}
\item[\textup{(ii)}] For all $h \in L^2(\Gamma)$, we have
\begin{equation}\label{eq:bessel}
\sum_{k \geq 1} \mu_k^{-1} \abs{\ps{h}{\psi_k}_\Gamma}^2
\leq \frn^2 \kappa^{-1} \norm{h}_{L^2(\Gamma)}^2.
\end{equation}
\end{enumerate}
\end{lemma}

\begin{proof}
(i) It follows from \eqref{eq:beig} that
$\frb(\phi_k,\phi_j) = \mu_k \delta_{jk}$ for all $j,\, k \in \NN$, where
$\delta_{jk}$ is the Kronecker symbol. Thus $\{ \mu_k^{-1/2} \phi_k,\ k \in \NN \}$
is an orthonormal system of the Hilbert space $(V,\frb)$. This system is
complete: if $v \in V$ satisfies $\frb(\phi_k,v) = 0$ for all $k \in \NN$, then
we have $(\phi_k,v) = 0$ for all $k \in \NN$ by \eqref{eq:beig}, and hence
$v = 0$. Therefore, $\{ \mu_k^{-1/2} \phi_k,\ k \in \NN \}$ is an orthonormal
basis of $(V,\frb)$. Now, let $u \in V$. Since $\frb$ is Hermitian and $\mu_k$
is real, we have
$\frb(u,\phi_k) = \overline{\frb(\phi_k,u)} = \mu_k (u,\phi_k)$ from
\eqref{eq:beig}, so the expansion of $u$ in this basis reads
\[
u = \sum_{k \geq 1} \frb(u,\mu_k^{-1/2} \phi_k)\, \mu_k^{-1/2} \phi_k
= \sum_{k \geq 1} (u,\phi_k) \phi_k,
\]
the series being convergent in $(V,\frb)$, hence in $V$ by \eqref{eq:b}, and
\eqref{eq:parseval} is nothing but the Parseval identity. Conversely, let
$u \in H$ be such that $\sum_{k \geq 1} \mu_k \abs{(u,\phi_k)}^2 < \infty$. Then
the series
$\sum_{k \geq 1} (u,\phi_k) \phi_k
= \sum_{k \geq 1} \mu_k^{1/2} (u,\phi_k)\, \mu_k^{-1/2} \phi_k$
converges in $V$ to some function $w \in V$. Since $V$ is continuously embedded
in $H$, the series converges to $w$ in $H$ as well. As it also converges to $u$
in $H$, the family $\{ \phi_k,\ k \in \NN \}$ being an orthonormal basis of $H$,
this yields $u = w \in V$.

(ii) Fix $N \in \NN$, put $c_k := \ps{h}{\psi_k}_\Gamma$ for $k \in \NN$, and set
$v_N := \sum_{k=1}^N \mu_k^{-1} c_k \phi_k \in V$. On the one hand, since the
$\psi_k$ are real-valued, we have
\[
\ps{h}{{v_N}_{|\Gamma}}_\Gamma
= \sum_{k=1}^N \mu_k^{-1} \overline{c_k}\, \ps{h}{\psi_k}_\Gamma
= \sum_{k=1}^N \mu_k^{-1} \abs{c_k}^2 =: S_N,
\]
and on the other hand, \eqref{eq:b} and \eqref{eq:beig} yield
\[
\kappa \norm{v_N}_V^2 \leq \frb(v_N,v_N)
= \sum_{k=1}^N \mu_k\, \mu_k^{-2} \abs{c_k}^2 = S_N.
\]
Therefore, we get from the Cauchy--Schwarz inequality and \eqref{eq:trace} that
\[
S_N \leq \norm{h}_{L^2(\Gamma)}\, \frn \norm{v_N}_V
\leq \frn\, \kappa^{-1/2} \norm{h}_{L^2(\Gamma)} S_N^{1/2},
\]
whence $S_N \leq \frn^2 \kappa^{-1} \norm{h}_{L^2(\Gamma)}^2$. Since $N$ is
arbitrary, this entails \eqref{eq:bessel}.
\end{proof}

\subsection{Unique continuation and linear independence of the boundary traces}

Given $\lambda \in \RR$, we denote by $\ker(A - \lambda)$ the set of functions
$\phi \in V$ satisfying
\begin{equation}\label{eq:eigweak}
\fra(\phi,v) = \lambda (\phi,v), \quad v \in V.
\end{equation}
The following unique continuation result is the Robin counterpart of
\cite[Corollary~3.1]{Ch24b}. Notice that only the Dirichlet trace of $\phi$ is
involved in its statement: roughly speaking, the Neumann trace of $\phi$ vanishes
wherever its Dirichlet trace does, due to the homogeneous Robin boundary
condition which is encoded in \eqref{eq:eigweak}.

\begin{lemma}\label{lem:ucp}
Let $\lambda \in \RR$ and let $\phi \in \ker(A - \lambda)$. If
$\phi_{|\Gamma} = 0$ a.e.\ on a nonempty open subset $\Gamma'$ of $\Gamma$, then
$\phi = 0$ in $\Omega$.
\end{lemma}

\begin{proof}
Pick $x_0 \in \Gamma'$ and choose $r_0 > 0$ so small that the open ball $B$
centered at $x_0$ with radius $r_0$ satisfies $B \cap \Gamma \subset \Gamma'$.
Then $\widehat{\Omega} := \Omega \cup B$ is a connected open subset of $\RR^n$
and, the domain $\Omega$ being Lipschitz, $B \setminus \overline{\Omega}$ is a
nonempty open subset of $\widehat{\Omega}$. We denote by $\widehat{\phi}$
(resp., $\widehat{q}$) the extension by zero of $\phi$ (resp., $q$) to
$\widehat{\Omega}$.

Let $\chi \in C_c^\infty(\widehat{\Omega})$. Since
$\supp \chi \cap \Gamma \subset B \cap \Gamma \subset \Gamma'$ and
$\phi_{|\Gamma} = 0$ a.e.\ on $\Gamma'$, we get upon integrating by parts in
$\Omega$ that
\[
\int_{\widehat{\Omega}} \widehat{\phi}\, \partial_i \chi\, \dd x
= \int_\Omega \phi\, \partial_i \chi\, \dd x
= - \int_\Omega \partial_i \phi\, \chi\, \dd x
+ \int_\Gamma \phi\, \chi\, \nu_i\, \dd\sigma
= - \int_\Omega \partial_i \phi\, \chi\, \dd x, \quad i = 1,\ldots,n,
\]
where $\nu_i$ is the $i$-th component of $\nu$. Therefore, we have
$\widehat{\phi} \in H^1(\widehat{\Omega})$, and $\nabla \widehat{\phi}$ is the
extension by zero of $\nabla \phi$ to $\widehat{\Omega}$. Similarly, the
boundary term $\int_\Gamma \alpha \phi \chi\, \dd\sigma$ vanishes, so we obtain
by taking $v = \overline{\chi}_{|\Omega} \in V$ in \eqref{eq:eigweak}, that
\[
\int_{\widehat{\Omega}} \nabla \widehat{\phi} \cdot \nabla \chi\, \dd x
+ \int_{\widehat{\Omega}} (\widehat{q} - \lambda) \widehat{\phi}\, \chi\, \dd x
= 0, \quad \chi \in C_c^\infty(\widehat{\Omega}).
\]
This means that
\begin{equation}\label{eq:ext-eq}
-\Delta \widehat{\phi} = (\lambda - \widehat{q}) \widehat{\phi}
\ \text{ in the distributional sense in } \widehat{\Omega}.
\end{equation}

Further, as $\phi \in V \subset L^{p^*}(\Omega)$, we have
$\widehat{\phi} \in L^{p^*}(\widehat{\Omega})$, and since
$\widehat{q} - \lambda \in L^{n/2}(\widehat{\Omega})$ and $1/p = 2/n + 1/p^*$,
the right-hand side of \eqref{eq:ext-eq} lies in $L^p(\widehat{\Omega})$ by
H\"older's inequality. Thus we have
$\widehat{\phi} \in W^{2,p}_{\mathrm{loc}}(\widehat{\Omega})$ by interior
elliptic regularity. Indeed, for any open ball $B'$ whose closure is contained
in $\widehat{\Omega}$, the Newtonian potential $w$ of the restriction of
$\Delta \widehat{\phi}$ to $B'$ belongs to $W^{2,p}(B')$ and satisfies
$\Delta w = \Delta \widehat{\phi}$ in $B'$, according to
\cite[Theorem~9.9]{GT01}, so $\widehat{\phi} - w$ is harmonic in $B'$ in the
distributional sense, and hence smooth in $B'$ by Weyl's lemma. Therefore,
$\widehat{\phi} \in W^{2,p}(B'')$ for every open ball $B''$ whose closure is
contained in $B'$, and the claim follows since every compact subset of
$\widehat{\Omega}$ can be covered by finitely many such balls $B''$.

Summing up, $\widehat{\phi} \in W^{2,p}_{\mathrm{loc}}(\widehat{\Omega})$
satisfies
$\abs{\Delta \widehat{\phi}} \leq \abs{\widehat{q} - \lambda}
\abs{\widehat{\phi}}$ a.e.\ in the connected open set $\widehat{\Omega}$, with
$\widehat{q} - \lambda \in L^{n/2}(\widehat{\Omega})$, and it vanishes in the
nonempty open set $B \setminus \overline{\Omega}$. Therefore, we have
$\widehat{\phi} = 0$ in $\widehat{\Omega}$ by the unique continuation property
established in \cite[Theorem~6.3]{JK85}, and consequently $\phi = 0$ in
$\Omega$.
\end{proof}

\begin{corollary}\label{cor:indep}
Let $\lambda$ be an eigenvalue of $A$ with multiplicity $d \in \NN$, and let
$\{ \phi^1,\ldots,\phi^d \}$ be a basis of $\ker(A - \lambda)$. Then, for every
nonempty open subset $\Gamma'$ of $\Gamma$, the functions
$\phi^1_{|\Gamma'},\ldots,\phi^d_{|\Gamma'}$ are linearly independent in
$L^2(\Gamma')$.
\end{corollary}

\begin{proof}
Let $(c_1,\ldots,c_d) \in \CC^d$ be such that
$\sum_{j=1}^d c_j \phi^j_{|\Gamma'} = 0$ in $L^2(\Gamma')$. Then the function
$\phi := \sum_{j=1}^d c_j \phi^j$, which lies in $\ker(A - \lambda)$, satisfies
$\phi_{|\Gamma} = 0$ a.e.\ on $\Gamma'$, so we have $\phi = 0$ in $\Omega$ by
Lemma~\ref{lem:ucp}. Since $\phi^1,\ldots,\phi^d$ are linearly independent in
$H$, this entails that $c_1 = \cdots = c_d = 0$.
\end{proof}

\subsection{Fractional integrals and Mittag-Leffler functions}

Recall the definitions \eqref{eq:kbeta}--\eqref{eq:RL}. In the convolution
products below, the functions defined on $(0,T)$ are extended by zero outside
$(0,T)$, while the kernels $k_\beta$, $\beta \in (0,\infty)$, and $K_\lambda$,
$\lambda \in \RR$, defined in \eqref{eq:kbeta} and \eqref{eq:Klambda} on
$(0,\infty)$, are extended by zero on $(-\infty,0]$. With this convention, the
convolution of two such functions $w_1$ and $w_2$ reads
$(w_1 * w_2)(t) = \int_0^t w_1(t-\tau)\, w_2(\tau)\, \dd\tau$ for $t > 0$, in
accordance with \eqref{eq:RL}.

\begin{lemma}\label{lem:RL}
Let $\beta,\, \beta' \in (0,\infty)$. Then the following statements hold true.
\begin{enumerate}
\item[\textup{(i)}] $k_\beta * k_{\beta'} = k_{\beta+\beta'}$ in $(0,\infty)$.
As a consequence, $I^\beta I^{\beta'} w = I^{\beta+\beta'} w$ for all
$w \in L^1(0,T)$.
\item[\textup{(ii)}] For all $K \in L^1(0,T)$ and $w \in L^2(0,T)$, we have
\[
\norm{K * w}_{L^2(0,T)} \leq \norm{K}_{L^1(0,T)} \norm{w}_{L^2(0,T)}.
\]
In particular, we have
\begin{equation}\label{eq:RLbound}
\norm{I^\beta w}_{L^2(0,T)} \leq \frac{T^\beta}{\Gamma(\beta+1)}
\norm{w}_{L^2(0,T)}, \quad w \in L^2(0,T).
\end{equation}
\end{enumerate}
\end{lemma}

\begin{proof}
(i) For $t > 0$, the change of variable $\tau = t s$ yields
\begin{align*}
(k_\beta * k_{\beta'})(t)
& = \frac{1}{\Gamma(\beta) \Gamma(\beta')} \int_0^t (t-\tau)^{\beta-1}
\tau^{\beta'-1}\, \dd\tau \\
& = \frac{t^{\beta+\beta'-1}}{\Gamma(\beta) \Gamma(\beta')}
\int_0^1 (1-s)^{\beta-1} s^{\beta'-1}\, \dd s
= k_{\beta+\beta'}(t),
\end{align*}
since the last integral is the Beta function
$B(\beta,\beta') = \Gamma(\beta) \Gamma(\beta') / \Gamma(\beta+\beta')$. The
second claim follows from this and the associativity of the convolution product
of integrable functions, which is a consequence of Fubini's theorem.

(ii) The first claim is Young's convolution inequality, applied to the
extensions by zero of $K$ and $w$, and \eqref{eq:RLbound} follows from it since
$\norm{k_\beta}_{L^1(0,T)} = T^\beta / (\beta \Gamma(\beta))
= T^\beta / \Gamma(\beta+1)$.
\end{proof}

For $\beta \in (0,\infty)$, the two-parameter Mittag-Leffler function is the
entire function defined by
\[
E_{\rho,\beta}(z) := \sum_{j \geq 0} \frac{z^j}{\Gamma(\rho j + \beta)},
\quad z \in \CC,
\]
see e.g.\ \cite[Section~1.2]{Pod99} or \cite{GKMR20}. Notice that $E_{1,1}(z) = e^z$.
For $\lambda \in \RR$ we set
\begin{equation}\label{eq:Klambda}
K_\lambda(t) := t^{\rho-1} E_{\rho,\rho}(-\lambda t^\rho)
= \sum_{j \geq 0} (-\lambda)^j k_{(j+1)\rho}(t), \quad t \in (0,\infty).
\end{equation}
For every $\lambda \in \RR$ and every $T' \in (0,\infty)$, we have
\begin{equation}\label{eq:KL1}
\sum_{j \geq 0} \abs{\lambda}^j \norm{k_{(j+1)\rho}}_{L^1(0,T')}
= \sum_{j \geq 0} \frac{\abs{\lambda}^j\, T'^{(j+1)\rho}}{\Gamma((j+1)\rho+1)}
= T'^\rho E_{\rho,\rho+1}(\abs{\lambda} T'^\rho) < \infty,
\end{equation}
so the series in \eqref{eq:Klambda} converges absolutely in $L^1(0,T')$, and
$K_\lambda \in L^1(0,T')$. The function $K_\lambda$ is the resolvent kernel of
the Volterra equation associated with \eqref{eq:weak}, in the following sense.

\begin{lemma}\label{lem:volterra}
Let $\lambda \in \RR$ and $\gamma \in L^2(0,T)$. Then $y := K_\lambda * \gamma$ is
the unique function $y \in L^2(0,T)$ satisfying
\begin{equation}\label{eq:volterra}
y + \lambda I^\rho y = I^\rho \gamma \quad \text{a.e.\ in } (0,T).
\end{equation}
\end{lemma}

\begin{proof}
\emph{Existence.} Since $K_\lambda \in L^1(0,T)$, the function
$y = K_\lambda * \gamma$ lies in $L^2(0,T)$ by Lemma~\ref{lem:RL}(ii). Next, the
series
in \eqref{eq:Klambda} being absolutely convergent in $L^1(0,T)$, we infer from
Lemma~\ref{lem:RL}(i) that
\[
\lambda\, k_\rho * K_\lambda
= \sum_{j \geq 0} (-\lambda)^j \lambda\, k_{(j+2)\rho}
= - \sum_{j \geq 1} (-\lambda)^j k_{(j+1)\rho}
= - \big( K_\lambda - k_\rho \big)
\quad \text{in } L^1(0,T),
\]
i.e.\ $K_\lambda + \lambda\, k_\rho * K_\lambda = k_\rho$. Convolving this
identity with $\gamma$ and using the associativity of the convolution product,
we obtain $y + \lambda\, k_\rho * y = k_\rho * \gamma$, which is
\eqref{eq:volterra}.

\emph{Uniqueness.} By linearity, it is enough to show that any function
$y \in L^2(0,T)$ satisfying $y = -\lambda I^\rho y$ vanishes identically.
Iterating this identity and using Lemma~\ref{lem:RL}(i), we get that
$y = (-\lambda)^N I^{N\rho} y$ for all $N \in \NN$. Hence we have
\[
\norm{y}_{L^2(0,T)} \leq \frac{\abs{\lambda}^N T^{N\rho}}{\Gamma(N\rho+1)}
\norm{y}_{L^2(0,T)}, \quad N \in \NN,
\]
by \eqref{eq:RLbound}. Since $\Gamma(N\rho+1)$ grows faster than any geometric
sequence as $N \to \infty$, the right-hand side tends to zero, and $y = 0$.
\end{proof}

We turn now to estimating $K_\lambda$. In the remaining part of this text, we
fix an angle
\[
\omega \in \big( \pi\rho/2,\, \min\{\pi,\pi\rho\} \big),
\quad \text{and we set} \quad
\vartheta := \min\{ \pi/2,\, (\pi - \omega)/\rho \} \in (0,\pi/2].
\]
We denote by $\cS := \{ z \in \CC \setminus \{0\};\ \abs{\arg z} < \vartheta \}$
the open sector of half-angle $\vartheta$ around the positive real axis, where
$\arg z \in (-\pi,\pi]$ is the principal argument of $z$. For $z \in \cS$, the
powers $z^\rho$ and $z^{\rho-1}$ are defined with the principal branch of the
logarithm, and we set
$K_\lambda(z) := z^{\rho-1} E_{\rho,\rho}(-\lambda z^\rho)$, which is consistent
with \eqref{eq:Klambda} for $z \in (0,\infty)$. As $E_{\rho,\rho}$ is entire and
the principal branch of the logarithm is holomorphic in $\cS$, the function
$z \mapsto K_\lambda(z)$ is holomorphic in $\cS$ for every $\lambda \in \RR$.

\begin{lemma}\label{lem:ML}
\begin{enumerate}
\item[\textup{(a)}] There exists a constant $C_0 > 0$, depending only on $\rho$
and $\omega$, such that
\begin{equation}\label{eq:MLbound}
\abs{E_{\rho,\rho}(z)} \leq \frac{C_0}{1 + \abs{z}^2},
\quad z \in \CC \setminus \{0\},\ \omega \leq \abs{\arg z} \leq \pi.
\end{equation}
\item[\textup{(b)}] Let $\lambda \in (0,\infty)$. Then we have
\begin{equation}\label{eq:Kpos}
\abs{K_\lambda(t)} \leq \frac{C_0\, t^{\rho-1}}{1 + \lambda^2 t^{2\rho}},
\quad t \in (0,\infty),
\end{equation}
and consequently
\begin{equation}\label{eq:Kpos2}
\int_0^\infty \abs{K_\lambda(t)}\, \dd t \leq \frac{\pi C_0}{2 \rho \lambda}
\quad \text{and} \quad
\abs{K_\lambda(t)} \leq \frac{C_0}{\lambda^2\, t^{\rho+1}},
\quad t \in (0,\infty).
\end{equation}
Moreover, we have
\begin{equation}\label{eq:Ksector}
\abs{K_\lambda(z)} \leq \frac{C_0}{2 \lambda \abs{z}}, \quad z \in \cS.
\end{equation}
\item[\textup{(c)}] Let $\Lambda \in (0,\infty)$, $T' \in (0,\infty)$ and
$\delta \in (0,T')$. Then there exist a constant $C_1 > 0$, depending only on
$\rho$, $\omega$, $\Lambda$ and $T'$, and a constant $C_2 > 0$, depending only on
$\rho$, $\omega$, $\Lambda$, $T'$ and $\delta$, such that for every
$\lambda \in (-\Lambda,\infty)$, we have
\begin{equation}\label{eq:Kunif}
\norm{K_\lambda}_{L^1(0,T')} \leq \frac{C_1}{\lambda + \Lambda}
\quad \text{and} \quad
\sup_{t \in [\delta,T']} \abs{K_\lambda(t)} \leq \frac{C_2}{(\lambda+\Lambda)^2}.
\end{equation}
\item[\textup{(d)}] Let $\lambda \in \RR$ and let $p \in (0,\infty)$ satisfy
$p^\rho + \lambda > 0$. Then $t \mapsto e^{-pt} K_\lambda(t)$ is integrable on
$(0,\infty)$, and we have
\begin{equation}\label{eq:laplace}
\int_0^\infty e^{-pt} K_\lambda(t)\, \dd t = \frac{1}{p^\rho + \lambda}.
\end{equation}
\end{enumerate}
\end{lemma}

\begin{proof}
(a) Since $\rho \in (0,2)$ and $\pi\rho/2 < \omega < \min\{\pi,\pi\rho\}$, the
asymptotic expansion \cite[Theorem~1.4, Eq.~(1.143)]{Pod99}, applied with
$\alpha = \beta = \rho$ and $p = 1$ in the notations of \cite{Pod99}, yields
\[
E_{\rho,\rho}(z) = - \frac{z^{-1}}{\Gamma(\rho - \rho)}
+ O\big( \abs{z}^{-2} \big)
\quad \text{as } \abs{z} \to \infty,\ \omega \leq \abs{\arg z} \leq \pi,
\]
where $1/\Gamma(0) := 0$, since $1/\Gamma$ is an entire function vanishing at
$0$ (see also \cite[Eq.~(2.6)]{KOSY18} for the same estimate). Hence there exist
two constants $R \in [1,\infty)$ and $C' > 0$ such that
$\abs{E_{\rho,\rho}(z)} \leq C' \abs{z}^{-2} \leq 2 C' / (1 + \abs{z}^2)$ for all
$z$ with $\abs{z} \geq R$ and $\omega \leq \abs{\arg z} \leq \pi$, where we used
that $1 + \abs{z}^2 \leq 2 \abs{z}^2$ for $\abs{z} \geq 1$. Moreover, putting
$C'' := \max_{\abs{z} \leq R} \abs{E_{\rho,\rho}(z)}$, we have
$\abs{E_{\rho,\rho}(z)} \leq C'' (1+R^2) / (1 + \abs{z}^2)$ for
$\abs{z} \leq R$. Thus \eqref{eq:MLbound} holds with
$C_0 := \max\{ 2C',\, (1+R^2) C'' \}$.

(b) For $t > 0$, the complex number $-\lambda t^\rho$ is negative, hence
$\arg(-\lambda t^\rho) = \pi$, and \eqref{eq:Kpos} follows from \eqref{eq:MLbound}.
Next, the change of variable $x = \lambda t^\rho$, for which
$\dd x = \rho \lambda t^{\rho-1} \dd t$, gives
\[
\int_0^\infty \frac{t^{\rho-1}}{1 + \lambda^2 t^{2\rho}}\, \dd t
= \frac{1}{\rho \lambda} \int_0^\infty \frac{\dd x}{1 + x^2}
= \frac{\pi}{2 \rho \lambda},
\]
which, together with \eqref{eq:Kpos}, yields the first estimate of
\eqref{eq:Kpos2}. The second one follows from \eqref{eq:Kpos} and the inequality
$1 + \lambda^2 t^{2\rho} \geq \lambda^2 t^{2\rho}$.

Let us now prove \eqref{eq:Ksector}. Let $z \in \cS$ and put
$\varphi := \arg z \in (-\vartheta,\vartheta)$. Then
$-\lambda z^\rho = \lambda \abs{z}^\rho e^{i(\rho\varphi + \pi)}$, where
$\rho \abs{\varphi} < \rho\vartheta \leq \pi - \omega$. Thus, if $\varphi > 0$,
the principal argument of $-\lambda z^\rho$ is $\rho\varphi - \pi
\in (-\pi,-\omega)$, and if $\varphi \leq 0$, it is $\rho\varphi + \pi
\in (\omega,\pi]$. In both cases we have
$\abs{\arg(-\lambda z^\rho)} = \pi - \rho\abs{\varphi} \in (\omega,\pi]$, so we
get from \eqref{eq:MLbound} that
\[
\abs{K_\lambda(z)} \leq \frac{C_0 \abs{z}^{\rho-1}}{1 + \lambda^2 \abs{z}^{2\rho}}
\leq \frac{C_0 \abs{z}^{\rho-1}}{2 \lambda \abs{z}^\rho}
= \frac{C_0}{2 \lambda \abs{z}},
\]
where we used that $1 + x^2 \geq 2x$ for all $x \geq 0$.

(c) We distinguish two cases. If $\lambda \geq \Lambda$, then
$\lambda + \Lambda \leq 2\lambda$, and \eqref{eq:Kpos2} yields
\[
\norm{K_\lambda}_{L^1(0,T')} \leq \frac{\pi C_0}{2\rho\lambda}
\leq \frac{\pi C_0}{\rho (\lambda + \Lambda)}
\quad \text{and} \quad
\sup_{t \in [\delta,T']} \abs{K_\lambda(t)} \leq \frac{C_0}{\lambda^2
\delta^{\rho+1}} \leq \frac{4 C_0}{\delta^{\rho+1} (\lambda+\Lambda)^2}.
\]
If $-\Lambda < \lambda < \Lambda$, then $\abs{\lambda} < \Lambda$ and
$\lambda + \Lambda < 2\Lambda$. We infer from \eqref{eq:Klambda} and
\eqref{eq:KL1} that
\[
\norm{K_\lambda}_{L^1(0,T')} \leq T'^\rho E_{\rho,\rho+1}(\Lambda T'^\rho)
\leq \frac{2\Lambda\, T'^\rho E_{\rho,\rho+1}(\Lambda T'^\rho)}{\lambda+\Lambda}.
\]
Similarly, since the coefficients of the power series defining $E_{\rho,\rho}$
are positive, we have for all $t \in [\delta,T']$,
\[
\abs{K_\lambda(t)} \leq t^{\rho-1} \sum_{j \geq 0}
\frac{\abs{\lambda}^j t^{\rho j}}{\Gamma(\rho j + \rho)}
= t^{\rho-1} E_{\rho,\rho}(\abs{\lambda} t^\rho)
\leq d_\rho\, E_{\rho,\rho}(\Lambda T'^\rho)
\leq \frac{4 \Lambda^2 d_\rho\, E_{\rho,\rho}(\Lambda T'^\rho)}
{(\lambda+\Lambda)^2},
\]
where $d_\rho := \max\{\delta^{\rho-1}, T'^{\rho-1}\}$. Therefore,
\eqref{eq:Kunif} holds with
\[
C_1 := \max\big\{ \pi C_0/\rho,\
2\Lambda\, T'^\rho E_{\rho,\rho+1}(\Lambda T'^\rho) \big\}
\ \text{ and } \
C_2 := \max\big\{ 4 C_0 \delta^{-\rho-1},\
4 \Lambda^2 d_\rho\, E_{\rho,\rho}(\Lambda T'^\rho) \big\}.
\]

(d) Recall that $\int_0^\infty e^{-pt} k_\beta(t)\, \dd t = p^{-\beta}$ for all
$\beta > 0$ and $p > 0$. Assume first that $p^\rho > \abs{\lambda}$. Then the
monotone convergence theorem yields
\[
\int_0^\infty e^{-pt} \sum_{j \geq 0} \abs{\lambda}^j k_{(j+1)\rho}(t)\, \dd t
= \sum_{j \geq 0} \abs{\lambda}^j p^{-(j+1)\rho}
= \frac{1}{p^\rho - \abs{\lambda}} < \infty.
\]
In light of \eqref{eq:Klambda}, this shows that $e^{-p\cdot} K_\lambda$ is
integrable on $(0,\infty)$, and that the series in \eqref{eq:Klambda} may be
integrated term by term against $e^{-pt}$, which yields
\[
\int_0^\infty e^{-pt} K_\lambda(t)\, \dd t
= \sum_{j \geq 0} (-\lambda)^j p^{-(j+1)\rho}
= \frac{p^{-\rho}}{1 + \lambda p^{-\rho}} = \frac{1}{p^\rho + \lambda}.
\]
This settles the case $\lambda \leq 0$, since the condition $p^\rho + \lambda > 0$
then reads $p^\rho > \abs{\lambda}$. Assume now that $\lambda > 0$. Then
$K_\lambda \in L^1(0,\infty)$ by \eqref{eq:Kpos2}, so the function
$F(\zeta) := \int_0^\infty e^{-\zeta t} K_\lambda(t)\, \dd t$ is well defined and
holomorphic in the half-plane $\{ \zeta \in \CC;\ \re \zeta > 0 \}$. Moreover,
for $\re \zeta > 0$, we have $\abs{\arg \zeta^\rho} < \rho\pi/2 < \pi$, where
$\zeta^\rho$ is defined with the principal branch of the logarithm, so that
$\zeta^\rho + \lambda \neq 0$, and the function
$\zeta \mapsto 1/(\zeta^\rho + \lambda)$ is holomorphic in the same half-plane.
As we have just seen that $F(\zeta) = 1/(\zeta^\rho + \lambda)$ for all real
$\zeta > \lambda^{1/\rho}$, the identity theorem for holomorphic functions entails
that this equality holds in the whole half-plane, and in particular for all
$\zeta = p \in (0,\infty)$.
\end{proof}

\section{Analysis of the direct problem}\label{sec:direct}

We use the notations of Section~\ref{sec:prelim}, and we write
\[
K_k := K_{\lambda_k}, \quad k \in \NN.
\]
We denote by $C_1$ (resp., $C_2$) the constant of Lemma~\ref{lem:ML}(c)
associated with $\Lambda = \lambda^*$ and $T' = T$ (resp., and some
$\delta \in (0,T)$, which will be specified when needed). Since
$\lambda_k > -\lambda^*$, we get from \eqref{eq:Kunif} that
\begin{equation}\label{eq:Kk}
\norm{K_k}_{L^1(0,T)} \leq C_1 \mu_k^{-1}
\quad \text{and} \quad
\sup_{t \in [\delta,T]} \abs{K_k(t)} \leq C_2 \mu_k^{-2}, \quad k \in \NN.
\end{equation}

\subsection{Existence, uniqueness and representation of the weak solution}

\begin{proposition}\label{prop:wp}
For every $g \in L^2(\Sigma)$, the IBVP \eqref{eq:ibvp} admits a unique weak
solution $u = u_g$ in the sense of Definition~\ref{def:weak}. Moreover, we have
\begin{equation}\label{eq:rep}
u = \sum_{k \geq 1} u_k \phi_k, \qquad
u_k := K_k * \gamma_k, \qquad
\gamma_k(\tau) := \ps{g(\cdot,\tau)}{\psi_k}_\Gamma, \quad \tau \in (0,T),
\end{equation}
where the series converges in $L^2(0,T;V)$, and the estimate
\begin{equation}\label{eq:energy}
\norm{u}_{L^2(0,T;V)} \leq C_1\, \frn\, \kappa^{-1} \norm{g}_{L^2(\Sigma)}
\end{equation}
holds.
\end{proposition}

\begin{proof}
For all $k \in \NN$, we have $\gamma_k \in L^2(0,T)$ by the Cauchy--Schwarz
inequality, and we infer from \eqref{eq:bessel} that
\begin{equation}\label{eq:gammasum}
\sum_{k \geq 1} \mu_k^{-1} \abs{\gamma_k(\tau)}^2
\leq \frn^2 \kappa^{-1} \norm{g(\cdot,\tau)}_{L^2(\Gamma)}^2
\quad \text{for a.e. } \tau \in (0,T).
\end{equation}

\emph{Step 1: Uniqueness and representation.} Let $u$ be a weak solution to
\eqref{eq:ibvp} and put $y_k(t) := (u(\cdot,t),\phi_k)$ for a.e.\
$t \in (0,T)$ and $k \in \NN$, in such a way that $y_k \in L^2(0,T)$. Since
$\fra$ is Hermitian and $\lambda_k$ is real, we have
$\fra(u(\cdot,t),\phi_k) = \overline{\fra(\phi_k,u(\cdot,t))}
= \lambda_k\, y_k(t)$ for a.e.\ $t \in (0,T)$, according to \eqref{eq:eig}.
Moreover, as $\psi_k$ is real-valued, we have
$\ps{g(\cdot,t)}{\psi_k}_\Gamma = \gamma_k(t)$. Thus, taking $v = \phi_k$ in
\eqref{eq:weak}, we find that
\[
y_k + \lambda_k I^\rho y_k = I^\rho \gamma_k \quad \text{a.e.\ in } (0,T),
\]
and consequently $y_k = K_k * \gamma_k = u_k$, by Lemma~\ref{lem:volterra}.
Since $\{ \phi_k,\ k \in \NN \}$ is an orthonormal basis of $H$, this shows that
$u(\cdot,t) = \sum_{k \geq 1} u_k(t) \phi_k$ in $H$ for a.e.\ $t \in (0,T)$. In
particular, there is at most one weak solution to \eqref{eq:ibvp}.

\emph{Step 2: Convergence.} Let $u_k$, $k \in \NN$, be defined by
\eqref{eq:rep}. We get from Lemma~\ref{lem:RL}(ii) and \eqref{eq:Kk} that
\begin{equation}\label{eq:ukL2}
\norm{u_k}_{L^2(0,T)} \leq \norm{K_k}_{L^1(0,T)} \norm{\gamma_k}_{L^2(0,T)}
\leq C_1 \mu_k^{-1} \norm{\gamma_k}_{L^2(0,T)}, \quad k \in \NN.
\end{equation}
For $N \in \NN$ and a.e.\ $\tau \in (0,T)$, we put
$R_N(\tau) := \sum_{k \geq N} \mu_k^{-1} \abs{\gamma_k(\tau)}^2$. By
\eqref{eq:gammasum}, we have
\[
R_N(\tau) \leq \frn^2 \kappa^{-1} \norm{g(\cdot,\tau)}_{L^2(\Gamma)}^2
\quad \text{and} \quad \lim_{N \to \infty} R_N(\tau) = 0
\quad \text{for a.e. } \tau \in (0,T),
\]
so we have $\lim_{N \to \infty} \int_0^T R_N(\tau)\, \dd\tau = 0$ by the
dominated convergence theorem. Next, we denote by
$S_N := \sum_{k=1}^N u_k \phi_k$ the $N$-th partial sum of the series in
\eqref{eq:rep}, which belongs to $L^2(0,T;V)$. For all integers
$N' > N \geq 1$, we deduce from \eqref{eq:b}, \eqref{eq:parseval} and
\eqref{eq:ukL2} that
\begin{align*}
\kappa \int_0^T \norm{S_{N'}(t) - S_N(t)}_V^2\, \dd t
& \leq \sum_{k=N+1}^{N'} \mu_k \norm{u_k}_{L^2(0,T)}^2 \\
& \leq C_1^2 \sum_{k=N+1}^{N'} \mu_k^{-1} \norm{\gamma_k}_{L^2(0,T)}^2
\leq C_1^2 \int_0^T R_{N+1}(\tau)\, \dd\tau.
\end{align*}
Therefore, $(S_N)_{N \in \NN}$ is a Cauchy sequence in $L^2(0,T;V)$, and we
denote by $u \in L^2(0,T;V)$ its limit. Moreover, arguing as above with $S_N$
instead of $S_{N'} - S_N$, and then sending $N$ to infinity, we get that
\[
\kappa \norm{u}_{L^2(0,T;V)}^2 \leq C_1^2 \int_0^T R_1(\tau)\, \dd\tau
\leq C_1^2\, \frn^2 \kappa^{-1} \norm{g}_{L^2(\Sigma)}^2,
\]
which yields \eqref{eq:energy}.

\emph{Step 3: Existence.} It remains to show that the function $u$ built in
Step 2 fulfills \eqref{eq:weak}. To this end, we fix $v \in V$, we put
$v_k := (v,\phi_k)$ for $k \in \NN$, and we set
$P_N v := \sum_{k=1}^N v_k \phi_k$ for $N \in \NN$. Since
$(\phi_k,v) = \overline{v_k}$ and $\fra(\phi_k,v) = \lambda_k \overline{v_k}$
by \eqref{eq:eig}, and since $u_k + \lambda_k I^\rho u_k = I^\rho \gamma_k$ a.e.\
in $(0,T)$ by Lemma~\ref{lem:volterra}, we have
\begin{equation}\label{eq:weakN}
\begin{aligned}
(S_N,v) + I^\rho\big[ \fra(S_N,v) \big]
& = \sum_{k=1}^N \overline{v_k} \big( u_k + \lambda_k I^\rho u_k \big)
= \sum_{k=1}^N \overline{v_k}\, I^\rho \gamma_k \\
& = I^\rho \big[ \ps{g}{(P_N v)_{|\Gamma}}_\Gamma \big]
\quad \text{a.e.\ in } (0,T),
\end{aligned}
\end{equation}
where we used in the last identity that
$\ps{g(\cdot,\tau)}{(P_N v)_{|\Gamma}}_\Gamma = \sum_{k=1}^N \overline{v_k}\,
\gamma_k(\tau)$, the $\psi_k$ being real-valued. Let us send $N$ to infinity in
\eqref{eq:weakN}, in the space $L^2(0,T)$. First, the Cauchy--Schwarz inequality
yields
$\norm{(S_N - u,v)}_{L^2(0,T)} \leq \norm{S_N - u}_{L^2(0,T;V)} \norm{v}_V$,
and we get from \eqref{eq:acont} and \eqref{eq:RLbound} that
\[
\norm{I^\rho\big[ \fra(S_N - u,v) \big]}_{L^2(0,T)}
\leq \frac{T^\rho}{\Gamma(\rho+1)}\, C_\fra \norm{S_N - u}_{L^2(0,T;V)}
\norm{v}_V.
\]
Both right-hand sides tend to zero as $N \to \infty$, by Step 2. Second, since
$P_N v$ converges to $v$ in $V$ by Lemma~\ref{lem:spec}(i), we infer from
\eqref{eq:trace} and \eqref{eq:RLbound} that
\[
\norm{I^\rho\big[ \ps{g}{(P_N v - v)_{|\Gamma}}_\Gamma \big]}_{L^2(0,T)}
\leq \frac{T^\rho}{\Gamma(\rho+1)}\, \frn \norm{g}_{L^2(\Sigma)}
\norm{P_N v - v}_V \underset{N \to \infty}{\longrightarrow} 0.
\]
Therefore, \eqref{eq:weakN} yields \eqref{eq:weak} in the limit
$N \to \infty$.
\end{proof}

\begin{remark}\label{rem:caputo}
Let us justify the derivation of \eqref{eq:weak} sketched in
Section~\ref{sec:RtD}. Let $\rho \in (0,1) \cup (1,2)$ and let
$y \in C^m([0,T])$ satisfy $y^{(j)}(0) = 0$ for $j = 0,\ldots,m-1$. Since
$\partial_t^\rho y = I^{m-\rho} y^{(m)}$ by \eqref{eq:caputo}, we get from
Lemma~\ref{lem:RL}(i) that
\[
I^\rho \partial_t^\rho y = I^m y^{(m)}
= y - \sum_{j=0}^{m-1} \frac{y^{(j)}(0)}{j!}\, t^j = y,
\]
where the second identity is the Taylor formula with integral remainder, since
$I^m w(t) = \frac{1}{(m-1)!} \int_0^t (t-\tau)^{m-1} w(\tau)\, \dd\tau$.
The same is obviously true for $\rho = 1$. As a consequence, if
$u \in C^m([0,T];V)$ satisfies $\partial_t^j u(\cdot,0) = 0$ for
$j = 0,\ldots,m-1$, and
\[
\partial_t^\rho (u(\cdot,t),v) + \fra(u(\cdot,t),v)
= \ps{g(\cdot,t)}{v_{|\Gamma}}_\Gamma, \quad t \in (0,T),\ v \in V,
\]
then, applying $I^\rho$ to both sides of this identity, we find that $u$ is the
weak solution to \eqref{eq:ibvp}.
\end{remark}

\subsection{The solution at time \texorpdfstring{$T_0$}{T0}}

The following result shows that the weak solution to \eqref{eq:ibvp} is more
regular near $T_0$ when the excitation is switched off shortly before $T_0$.

\begin{proposition}\label{prop:T0}
Let $g \in L^2(\Sigma)$ vanish a.e.\ in
$\Gamma \times (T_0 - \varepsilon,T_0)$ for some $\varepsilon \in (0,T_0)$, let
$u = u_g$ be the weak solution to \eqref{eq:ibvp}, and put
$J_\varepsilon := [T_0 - \varepsilon/2,\, T_0]$. For $k \in \NN$ and
$t \in J_\varepsilon$, set
\begin{equation}\label{eq:ukT0}
U_k(t) := \int_0^{T_0 - \varepsilon} K_k(t-\tau)\, \gamma_k(\tau)\, \dd\tau,
\end{equation}
where $\gamma_k$ is defined in \eqref{eq:rep}. Then the series
$\sum_{k \geq 1} U_k \phi_k$ converges in $C(J_\varepsilon;V)$, and its sum
coincides with $u$ a.e.\ in $J_\varepsilon$. In particular, the restriction of
$u$ to $J_\varepsilon$ has a unique representative in $C(J_\varepsilon;V)$, its
value $u(\cdot,T_0)$ at time $T_0$ belongs to $V$, and we have
\begin{equation}\label{eq:traceT0}
u(\cdot,T_0)_{|\Gamma} = \sum_{k \geq 1} U_k(T_0)\, \psi_k
\quad \text{in } L^2(\Gamma).
\end{equation}
\end{proposition}

\begin{proof}
We use \eqref{eq:Kk} with $\delta := \varepsilon/2 \in (0,T)$. For
$t \in J_\varepsilon$ and $\tau \in (0,T_0-\varepsilon)$, we have
$t - \tau \in [\varepsilon/2, T]$, so the function
$(t,\tau) \mapsto K_k(t - \tau)$ is continuous and bounded by
$C_2 \mu_k^{-2}$ on $J_\varepsilon \times (0,T_0-\varepsilon)$. Since
$\gamma_k \in L^1(0,T)$, it follows from the dominated convergence theorem that
$U_k \in C(J_\varepsilon)$, and the Cauchy--Schwarz inequality yields
\begin{equation}\label{eq:Ukbound}
\sup_{t \in J_\varepsilon} \abs{U_k(t)} \leq C_2 \mu_k^{-2}\, T^{1/2}
\norm{\gamma_k}_{L^2(0,T)}, \quad k \in \NN.
\end{equation}
Next, since $\gamma_k$ vanishes a.e.\ in $(T_0-\varepsilon,T_0)$, we have
$(K_k * \gamma_k)(t) = U_k(t)$ for all $t \in J_\varepsilon$, so that
$u_k = U_k$ a.e.\ in $J_\varepsilon$, by \eqref{eq:rep}.

Further, for all integers $N' > N \geq 1$, we get from \eqref{eq:b},
\eqref{eq:parseval} and \eqref{eq:Ukbound} that
\[
\kappa \sup_{t \in J_\varepsilon}
\Big\lVert \sum_{k=N+1}^{N'} U_k(t) \phi_k \Big\rVert_V^2
\leq \sum_{k=N+1}^{N'} \mu_k \sup_{t \in J_\varepsilon} \abs{U_k(t)}^2
\leq C_2^2\, T \sum_{k \geq N+1} \mu_k^{-3} \norm{\gamma_k}_{L^2(0,T)}^2.
\]
As $\mu_k \geq \mu_1 > 0$, the right-hand side is bounded by
$C_2^2\, T \mu_1^{-2} \int_0^T R_{N+1}(\tau)\, \dd\tau$, where $R_N$ is defined
in the proof of Proposition~\ref{prop:wp}, and it tends to zero as
$N \to \infty$. Therefore, the series $\sum_{k \geq 1} U_k \phi_k$ converges in
$C(J_\varepsilon;V)$. Since the series $\sum_{k \geq 1} u_k \phi_k$ converges to
$u$ in $L^2(0,T;V)$ by Proposition~\ref{prop:wp}, and since $u_k = U_k$ a.e.\ in
$J_\varepsilon$ for every $k \in \NN$, the two sums coincide a.e.\ in
$J_\varepsilon$. A continuous representative being unique, the remaining claims
follow from this and the continuity of the trace operator from $V$ into
$L^2(\Gamma)$.
\end{proof}

As already mentioned in Section~\ref{sec:RtD}, every $g \in \cX_0$ fulfills the
assumption of Proposition~\ref{prop:T0} for some
$\varepsilon = \varepsilon(g) \in (0,T_0)$. Moreover, if this assumption is
satisfied by two numbers $\varepsilon' < \varepsilon$, the corresponding
continuous representatives coincide on $J_{\varepsilon'} \subset J_\varepsilon$,
so $u_g(\cdot,T_0)$ does not depend on the choice of $\varepsilon$. Therefore,
the Robin-to-Dirichlet map \eqref{eq:Dq} is well defined, and it
reads
\begin{equation}\label{eq:Dqrep}
D_q(g) = \sum_{k \geq 1} \Big( \int_0^{T_0-\varepsilon} K_k(T_0-\tau)
\ps{g(\cdot,\tau)}{\psi_k}_\Gamma\, \dd\tau \Big)\, {\psi_k}_{|\Gamma_1}
\ \text{ in } L^2(\Gamma_1), \quad g \in \cX_0.
\end{equation}

\section{Proof of Theorems \ref{thm:uniq} and \ref{thm:bsd}}\label{sec:proof}

Throughout this section, $\rho \in (0,2)$, $\Gamma_0$ and $\Gamma_1$ are two
nonempty open subsets of $\Gamma$ such that $\Gamma_0 \cap \Gamma_1 \neq
\emptyset$, and $q$, $\tq$ are two potentials in $L^{n/2}(\Omega,\RR)$. We denote
by $\lambda^*$ and $\kappa$ (resp., $\tla^*$ and $\widetilde{\kappa}$) the
constants of \eqref{eq:coer} associated with $q$ (resp., $\tq$), and we put
$\tK_k := K_{\tla_k}$ and $\widetilde{\mu}_k := \tla_k + \tla^*$, $k \in \NN$.
We also set
\[
\lambda_\sharp := \max\{ \lambda^*, \tla^* \}.
\]

We identify $L^2(\Gamma_0)$ with the subspace of $L^2(\Gamma)$ made of functions
vanishing a.e.\ in $\Gamma \setminus \Gamma_0$. For $f \in L^2(\Gamma)$ with
$\supp f \subset \Gamma_0$ and $h \in C_c^\infty(0,T_0)$, we denote by
$f \otimes h$ the function $(\sigma,t) \mapsto f(\sigma) h(t)$, extended by zero
to $\Gamma \times [T_0,T)$ when $T_0 < T$. It is clear that
$f \otimes h \in \cX_0$.

In view of Remark~\ref{rem:main}(ii), we consider a set $\cF$ of functions
$f \in L^2(\Gamma)$ with $\supp f \subset \Gamma_0$, which is dense in
$L^2(\Gamma_0)$, as well as a nonempty open interval $I \subset (0,T_0)$, and we
assume that
\begin{equation}\label{eq:hyp}
D_q(f \otimes h) = D_{\tq}(f \otimes h), \quad f \in \cF,\ h \in C_c^\infty(I).
\end{equation}
This condition is fulfilled under the assumption of Theorem~\ref{thm:bsd}, with
$I := (0,T_0)$ and
\[
\cF := \{ f \in L^2(\Gamma);\ \supp f \subset \Gamma_0 \}.
\]
This set is indeed dense in $L^2(\Gamma_0)$: as $\Gamma_0$ is the union of the
non-decreasing sequence of compact sets
\[
Z_j := \{ \sigma \in \Gamma;\ \abs{\sigma - \sigma'} \geq 1/j \text{ for all }
\sigma' \in \Gamma \setminus \Gamma_0 \}, \quad j \in \NN,
\]
each $f \in L^2(\Gamma_0)$ is the limit in $L^2(\Gamma_0)$ of the sequence
$(f \mathbf{1}_{Z_j})_{j \in \NN}$, where $\mathbf{1}_{Z_j}$ denotes the
characteristic function of $Z_j$, by the dominated convergence theorem.

The proof of Theorem~\ref{thm:bsd} is divided into five steps.

\subsection{Step 1: scalar generating functions}

For $f \in L^2(\Gamma)$ with $\supp f \subset \Gamma_0$ and
$\chi \in L^2(\Gamma_1)$, we set
\begin{equation}\label{eq:ck}
c_k := \ps{f}{\psi_k}_\Gamma\, \ps{\psi_k}{\chi}_{\Gamma_1}
= \int_{\Gamma_0 \times \Gamma_1} \psi_k(\sigma) \psi_k(\sigma')\,
f(\sigma)\, \overline{\chi(\sigma')}\, \dd\sigma\, \dd\sigma',
\quad k \in \NN,
\end{equation}
where the second identity follows from Fubini's theorem, since the $\psi_k$ are
real-valued and $f$ vanishes a.e.\ outside $\Gamma_0$. Let $\chi^0 \in
L^2(\Gamma)$ denote the extension of $\chi$ by zero to $\Gamma$. Since
$\abs{\ps{\psi_k}{\chi}_{\Gamma_1}} = \abs{\ps{\chi^0}{\psi_k}_\Gamma}$, we
infer from the Cauchy--Schwarz inequality and \eqref{eq:bessel}, applied with
$h = f$ and $h = \chi^0$, that
\begin{equation}\label{eq:cksum}
\sum_{k \geq 1} \mu_k^{-1} \abs{c_k}
\leq \Big( \sum_{k \geq 1} \mu_k^{-1} \abs{\ps{f}{\psi_k}_\Gamma}^2 \Big)^{1/2}
\Big( \sum_{k \geq 1} \mu_k^{-1} \abs{\ps{\chi^0}{\psi_k}_\Gamma}^2
\Big)^{1/2}
\leq \frn^2 \kappa^{-1} \norm{f}_{L^2(\Gamma)} \norm{\chi}_{L^2(\Gamma_1)}.
\end{equation}

Let $\NN_- := \{ k \in \NN;\ \lambda_k < \lambda^* \}$, which is a finite set
since $\lim_{k \to \infty} \lambda_k = +\infty$, and put
$\NN_+ := \NN \setminus \NN_-$. Notice that
\begin{equation}\label{eq:Nplus}
\mu_k = \lambda_k + \lambda^* \leq 2 \lambda_k, \quad k \in \NN_+.
\end{equation}

\begin{lemma}\label{lem:G}
Let $f \in L^2(\Gamma)$ with $\supp f \subset \Gamma_0$, let
$\chi \in L^2(\Gamma_1)$, and let $c_k$, $k \in \NN$, be defined by
\eqref{eq:ck}. Then the following statements hold true.
\begin{enumerate}
\item[\textup{(i)}] The series
\begin{equation}\label{eq:G}
G(z) := \sum_{k \geq 1} c_k\, K_k(z)
\end{equation}
converges absolutely and uniformly on every compact subset of the sector $\cS$,
and $G$ is holomorphic in $\cS$. In particular, $G$ is continuous on
$(0,\infty)$.
\item[\textup{(ii)}] For all $h \in C_c^\infty(0,T_0)$, we have
\begin{equation}\label{eq:DqG}
\ps{D_q(f \otimes h)}{\chi}_{\Gamma_1}
= \int_0^{T_0} h(T_0 - t)\, G(t)\, \dd t.
\end{equation}
\item[\textup{(iii)}] For all $p \in (0,\infty)$ such that
$p^\rho > \lambda^*$, the function $t \mapsto e^{-pt} G(t)$ is integrable on
$(0,\infty)$, the series $\sum_{k \geq 1} c_k / (p^\rho + \lambda_k)$ converges
absolutely, and we have
\begin{equation}\label{eq:LG}
\int_0^\infty e^{-pt}\, G(t)\, \dd t
= \sum_{k \geq 1} \frac{c_k}{p^\rho + \lambda_k}.
\end{equation}
\end{enumerate}
\end{lemma}

\begin{proof}
(i) Let $\cK$ be a compact subset of $\cS$, and put
$r_\cK := \min_{z \in \cK} \abs{z} > 0$. For $k \in \NN_+$, we have
$\lambda_k \geq \lambda^* > 0$, and we get from \eqref{eq:Ksector} and
\eqref{eq:Nplus} that
\[
\abs{c_k K_k(z)} \leq \frac{C_0 \abs{c_k}}{2 \lambda_k \abs{z}}
\leq \frac{C_0}{r_\cK}\, \mu_k^{-1} \abs{c_k}, \quad z \in \cK.
\]
In light of \eqref{eq:cksum}, the series $\sum_{k \in \NN_+} c_k K_k$ is thus
normally convergent on $\cK$. As $\NN_-$ is finite and each $K_k$ is holomorphic
in $\cS$, the series \eqref{eq:G} converges absolutely and uniformly on $\cK$,
and $G$ is holomorphic in $\cS$ by the Weierstrass theorem.

(ii) Let $h \in C_c^\infty(0,T_0)$, put $g := f \otimes h$, and pick
$\varepsilon \in (0,T_0/2)$ such that $\supp h \subset [\varepsilon,
T_0 - \varepsilon]$. Then $g$ vanishes a.e.\ in
$\Gamma \times (T_0-\varepsilon,T_0)$, and
$\ps{g(\cdot,\tau)}{\psi_k}_\Gamma = h(\tau) \ps{f}{\psi_k}_\Gamma$ for all
$\tau \in (0,T)$, where $h$ is extended by zero to $[T_0,T)$. Since $h$
vanishes in $(T_0 - \varepsilon,T_0)$, we infer from \eqref{eq:Dqrep}, upon
performing the change of variable $t = T_0 - \tau$, that
\begin{align*}
\ps{D_q(g)}{\chi}_{\Gamma_1}
& = \sum_{k \geq 1} \Big( \int_0^{T_0} K_k(T_0-\tau)\, h(\tau)\, \dd\tau \Big)
\ps{f}{\psi_k}_\Gamma \ps{\psi_k}{\chi}_{\Gamma_1} \\
& = \sum_{k \geq 1} c_k \int_0^{T_0} h(T_0-t)\, K_k(t)\, \dd t,
\end{align*}
where we used that the series in \eqref{eq:Dqrep} converges in
$L^2(\Gamma_1)$. Moreover, it follows from \eqref{eq:Kk} (with $T_0 \leq T$)
and \eqref{eq:cksum} that
\[
\sum_{k \geq 1} \abs{c_k} \int_0^{T_0} \abs{h(T_0-t)}\, \abs{K_k(t)}\, \dd t
\leq \norm{h}_{L^\infty(0,T_0)}\, C_1 \sum_{k \geq 1} \mu_k^{-1} \abs{c_k}
< \infty.
\]
Hence the sum over $k$ and the integral over $(0,T_0)$ may be interchanged, and
we obtain \eqref{eq:DqG}.

(iii) Let $p > 0$ satisfy $p^\rho > \lambda^*$. For $k \in \NN_+$, we get from
\eqref{eq:Kpos2} and \eqref{eq:Nplus} that
\[
\int_0^\infty e^{-pt} \abs{K_k(t)}\, \dd t
\leq \int_0^\infty \abs{K_k(t)}\, \dd t \leq \frac{\pi C_0}{2\rho\lambda_k}
\leq \frac{\pi C_0}{\rho}\, \mu_k^{-1}.
\]
For $k \in \NN_-$, we have $p^\rho + \lambda_k > \lambda^* + \lambda_k = \mu_k
> 0$, so $e^{-p\cdot} K_k$ is integrable on $(0,\infty)$ by
Lemma~\ref{lem:ML}(d). Since $\NN_-$ is finite, we deduce from this and
\eqref{eq:cksum} that
\[
\sum_{k \geq 1} \abs{c_k} \int_0^\infty e^{-pt} \abs{K_k(t)}\, \dd t < \infty.
\]
As a consequence, $e^{-p\cdot} G$ is integrable on $(0,\infty)$, and the series
\eqref{eq:G} may be integrated term by term against $e^{-pt}$. Since
$p^\rho + \lambda_k > 0$ for all $k \in \NN$, \eqref{eq:LG} follows from this and
\eqref{eq:laplace}. Finally, the absolute convergence of the series in
\eqref{eq:LG} follows from \eqref{eq:cksum}, since
$p^\rho + \lambda_k = (p^\rho - \lambda^*) + \mu_k \geq \mu_k$ for all
$k \in \NN$.
\end{proof}

\subsection{Step 2: analytic continuation}

Fix $f \in \cF$ and $\chi \in L^2(\Gamma_1)$. Let $c_k$ and $G$ be defined by
\eqref{eq:ck} and \eqref{eq:G}, and let $\tc_k$ and $\tG$ be defined in the same
way, upon substituting $(\tla_k,\tpsi_k)$ for $(\lambda_k,\psi_k)$. We shall
prove that
\begin{equation}\label{eq:GG}
G(t) = \tG(t), \quad t \in (0,\infty).
\end{equation}
Indeed, it follows from \eqref{eq:hyp} and from Lemma~\ref{lem:G}(ii), applied
to $q$ and to $\tq$, that
\[
\int_0^{T_0} h(T_0 - t)\, \big( G(t) - \tG(t) \big)\, \dd t = 0,
\quad h \in C_c^\infty(I).
\]
Put $I' := \{ T_0 - t;\ t \in I \}$, which is a nonempty open interval contained
in $(0,T_0)$. Since $t \mapsto h(T_0 - t)$ is arbitrary in $C_c^\infty(I')$ and
$G - \tG$ is continuous on $(0,\infty)$, this entails that $G(t) = \tG(t)$ for
all $t \in I'$. Further, $G - \tG$ is holomorphic in the connected open set
$\cS$, by Lemma~\ref{lem:G}(i), and it vanishes on the segment $I'$, which
has accumulation points in $\cS$. Therefore, $G - \tG$ vanishes identically in
$\cS$ by the identity theorem, which yields \eqref{eq:GG}.

\subsection{Step 3: Laplace transform and meromorphic functions}

Let $(\lambda_j')_{j \in \NN}$ be the strictly increasing sequence of the
distinct eigenvalues of $A$. For each $j \in \NN$, we denote by $m_j \in \NN$
the multiplicity of $\lambda_j'$, and by $\{ \phi_{j,l},\ l = 1,\ldots,m_j \}$
the family of those eigenfunctions $\phi_i$, $i \in \NN$, that are associated
with $\lambda_j'$, which is a real-valued orthonormal basis of
$\ker(A - \lambda_j')$. We set
\[
\Psi_j := (\psi_{j,1},\ldots,\psi_{j,m_j})^T,
\ \text{ where } \psi_{j,l} := {\phi_{j,l}}_{|\Gamma},\ l = 1,\ldots,m_j,
\]
and
\begin{equation}\label{eq:theta}
\theta_j(\sigma,\sigma') := \Psi_j(\sigma) \cdot \Psi_j(\sigma')
= \sum_{l=1}^{m_j} \psi_{j,l}(\sigma) \psi_{j,l}(\sigma'),
\quad (\sigma,\sigma') \in \Gamma \times \Gamma,
\end{equation}
where the symbol $\cdot$ stands for the Euclidean scalar product in
$\RR^{m_j}$. Notice that $\theta_j \in L^2(\Gamma \times \Gamma)$. Similarly, we
define $\tla_j'$, $\tm_j$, $\tphi_{j,l}$, $\tPsi_j$ and $\ttheta_j$ by replacing
$q$ by $\tq$ in the above definitions. Finally, we introduce the closed and
discrete subsets of the real line
\[
\cP := \{ -\lambda_j',\ j \in \NN \}
\quad \text{and} \quad
\widetilde{\cP} := \{ -\tla_j',\ j \in \NN \}.
\]

\begin{lemma}\label{lem:mero}
Let $f \in L^2(\Gamma)$ with $\supp f \subset \Gamma_0$, let
$\chi \in L^2(\Gamma_1)$, and let $c_k$, $k \in \NN$, be defined by
\eqref{eq:ck}. For $j \in \NN$, put
\begin{equation}\label{eq:Cj}
c_j' := \sum_{l=1}^{m_j} \ps{f}{\psi_{j,l}}_\Gamma\,
\ps{\psi_{j,l}}{\chi}_{\Gamma_1}
= \int_{\Gamma_0 \times \Gamma_1} \theta_j(\sigma,\sigma')\, f(\sigma)\,
\overline{\chi(\sigma')}\, \dd\sigma\, \dd\sigma'.
\end{equation}
Then the following statements hold true.
\begin{enumerate}
\item[\textup{(i)}] The series
\begin{equation}\label{eq:M}
M(s) := \sum_{k \geq 1} \frac{c_k}{s + \lambda_k}
= \sum_{j \geq 1} \frac{c_j'}{s + \lambda_j'}
\end{equation}
converge absolutely and uniformly on every compact subset of
$\CC \setminus \cP$, and $M$ is holomorphic in $\CC \setminus \cP$.
\item[\textup{(ii)}] For every $j \in \NN$, the function
$s \mapsto M(s) - c_j' / (s + \lambda_j')$ extends to a holomorphic function in a
neighborhood of $-\lambda_j'$.
\end{enumerate}
\end{lemma}

\begin{proof}
(i) Let $\cK$ be a compact subset of $\CC \setminus \cP$ and put
$R := \max_{s \in \cK} \abs{s}$. For all $k \in \NN$ such that
$\lambda_k \geq 2R + \lambda^*$, we have
\[
\abs{s + \lambda_k} \geq \lambda_k - R \geq \frac{\lambda_k + \lambda^*}{2}
= \frac{\mu_k}{2}, \quad s \in \cK,
\]
so $\abs{c_k / (s + \lambda_k)} \leq 2 \mu_k^{-1} \abs{c_k}$. As there are only
finitely many $k \in \NN$ such that $\lambda_k < 2R + \lambda^*$, and as each
function $s \mapsto c_k/(s + \lambda_k)$ is holomorphic in $\CC \setminus \cP$,
we deduce from \eqref{eq:cksum} that the first series in \eqref{eq:M} converges
absolutely and uniformly on $\cK$. The second series in \eqref{eq:M} is obtained
from the first one by grouping the terms associated with the same eigenvalue,
which is legitimate for absolutely convergent series. Indeed, by definition of
the $\phi_{j,l}$, we have
$\sum_{k \in \NN;\, \lambda_k = \lambda_j'} c_k = c_j'$ for all $j \in \NN$, and
the identity in \eqref{eq:Cj} follows from \eqref{eq:ck} and \eqref{eq:theta}.
Moreover, we have
$\abs{c_j' / (s + \lambda_j')} \leq \sum_{k;\, \lambda_k = \lambda_j'}
\abs{c_k / (s + \lambda_k)}$, so the second series also converges absolutely and
uniformly on $\cK$. Finally, $M$ is holomorphic in $\CC \setminus \cP$ by the
Weierstrass theorem.

(ii) Fix $j \in \NN$. By the same argument as in (i), the series
$\sum_{j' \neq j} c_{j'}' / (s + \lambda_{j'}')$ converges uniformly on every
compact subset of $\CC \setminus (\cP \setminus \{-\lambda_j'\})$, which is an
open neighborhood of $-\lambda_j'$. Its sum is thus holomorphic in this
neighborhood, and it coincides with $M(s) - c_j'/(s + \lambda_j')$ for
$s \notin \cP$.
\end{proof}

Fix $f \in \cF$ and $\chi \in L^2(\Gamma_1)$, and let $M$ and $c_j'$, $j \in \NN$,
be defined as in Lemma~\ref{lem:mero}. Let $\tM$ and $\tc_j'$, $j \in \NN$, be
defined in the same way upon substituting $\tq$ for $q$. For all $p > 0$ such that
$p^\rho > \lambda_\sharp$, we infer from \eqref{eq:GG} and from
Lemma~\ref{lem:G}(iii), applied to $q$ and to $\tq$, that
$M(p^\rho) = \tM(p^\rho)$. As $p^\rho$ runs over $(\lambda_\sharp,\infty)$ when
$p$ runs over $(\lambda_\sharp^{1/\rho},\infty)$, this means that
\begin{equation}\label{eq:MM0}
M(s) = \tM(s), \quad s \in (\lambda_\sharp,\infty).
\end{equation}
Notice that $(\lambda_\sharp,\infty) \cap (\cP \cup \widetilde{\cP}) =
\emptyset$, since $-\lambda_j' < \lambda^*$ and $-\tla_j' < \tla^*$ for all
$j \in \NN$. Now, the set $\cP \cup \widetilde{\cP}$ is a closed and discrete
subset of $\RR$, as $\lim_{j \to \infty} \lambda_j' = \lim_{j \to \infty}
\tla_j' = +\infty$, so $\CC \setminus (\cP \cup \widetilde{\cP})$ is a connected
open subset of $\CC$, and $M - \tM$ is holomorphic in it by
Lemma~\ref{lem:mero}(i). By the identity theorem, we deduce from \eqref{eq:MM0}
that
\begin{equation}\label{eq:MM}
M(s) = \tM(s), \quad s \in \CC \setminus (\cP \cup \widetilde{\cP}).
\end{equation}

\subsection{Step 4: identification of the eigenvalues and of the
\texorpdfstring{$\theta_j$}{theta j}}

We start by comparing the singularities of $M$ and $\tM$.

\begin{lemma}\label{lem:residue}
Let $f \in \cF$, $\chi \in L^2(\Gamma_1)$, and let $c_j'$ and $\tc_j'$,
$j \in \NN$, be as above. Let $j \in \NN$.
\begin{enumerate}
\item[\textup{(i)}] If $-\lambda_j' \notin \widetilde{\cP}$, then $c_j' = 0$.
\item[\textup{(ii)}] If $-\lambda_j' = -\tla_{l}'$ for some $l \in \NN$, then
$c_j' = \tc_{l}'$.
\end{enumerate}
\end{lemma}

\begin{proof}
Let $D$ be an open disc centered at $-\lambda_j'$ such that
$D \cap (\cP \cup \widetilde{\cP}) \subset \{-\lambda_j'\}$, and such that the
function $H_j(s) := M(s) - c_j'/(s + \lambda_j')$ extends holomorphically to $D$
(see Lemma~\ref{lem:mero}(ii)). We still denote by $H_j$ this extension.

(i) If $-\lambda_j' \notin \widetilde{\cP}$, then $\tM$ is holomorphic in $D$ by
Lemma~\ref{lem:mero}(i), and \eqref{eq:MM} yields
$c_j' / (s + \lambda_j') = \tM(s) - H_j(s)$ for all $s \in D \setminus
\{-\lambda_j'\}$. The right-hand side being bounded in a neighborhood of
$-\lambda_j'$, we necessarily have $c_j' = 0$.

(ii) If $-\lambda_j' = -\tla_{l}'$, then by Lemma~\ref{lem:mero}(ii) applied to
$\tM$, and upon possibly shrinking $D$, the function
$\tH_{l}(s) := \tM(s) - \tc_{l}'/(s + \tla_{l}')$ extends holomorphically to $D$.
Thus we get from \eqref{eq:MM} that
$(c_j' - \tc_{l}') / (s + \lambda_j') = \tH_{l}(s) - H_j(s)$ for all
$s \in D \setminus \{-\lambda_j'\}$, and hence $c_j' = \tc_{l}'$, for the same
reason as above.
\end{proof}

We next check that $\theta_j$ is not identically zero in
$\Gamma_0 \times \Gamma_1$, for all $j \in \NN$. We argue by contradiction and
assume that $\theta_j = 0$ a.e.\ in $\Gamma_0 \times \Gamma_1$. Then, for a.e.\
$\sigma \in \Gamma_0$, we have $\Psi_j(\sigma) \cdot \Psi_j = 0$ a.e.\ in
$\Gamma_1$ by Fubini's theorem, and hence $\Psi_j(\sigma) = 0$, since the
functions $\psi_{j,1},\ldots,\psi_{j,m_j}$ are linearly independent in
$L^2(\Gamma_1)$ according to Corollary~\ref{cor:indep}. Thus we have
$\psi_{j,1} = 0$ a.e.\ in $\Gamma_0$, which contradicts
Corollary~\ref{cor:indep} with $\Gamma' = \Gamma_0$. The same being true for
$\ttheta_j$, we have
\begin{equation}\label{eq:thetanz}
\theta_j \neq 0 \ \text{ and } \ \ttheta_j \neq 0
\ \text{ in } L^2(\Gamma_0 \times \Gamma_1), \quad j \in \NN.
\end{equation}

Further, the identity in \eqref{eq:Cj} can be rewritten as
\begin{equation}\label{eq:Cjps}
c_j' = \ps{\theta_j}{\overline{f} \otimes \chi}_{L^2(\Gamma_0 \times \Gamma_1)},
\end{equation}
where $\overline{f} \otimes \chi$ denotes the function
$(\sigma,\sigma') \mapsto \overline{f(\sigma)}\, \chi(\sigma')$. Since $\cF$ is
dense in $L^2(\Gamma_0)$, so is $\{ \overline{f},\ f \in \cF \}$. Moreover,
the linear span of the functions $a \otimes b : (\sigma,\sigma') \mapsto
a(\sigma) b(\sigma')$, where $a \in L^2(\Gamma_0)$ and $b \in L^2(\Gamma_1)$, is
dense in $L^2(\Gamma_0 \times \Gamma_1)$, and we have
$\norm{a \otimes b - \overline{f} \otimes b}_{L^2(\Gamma_0 \times \Gamma_1)} =
\norm{a - \overline{f}}_{L^2(\Gamma_0)} \norm{b}_{L^2(\Gamma_1)}$. Therefore, the
linear span of the functions $\overline{f} \otimes \chi$, where $f \in \cF$ and
$\chi \in L^2(\Gamma_1)$, is dense in $L^2(\Gamma_0 \times \Gamma_1)$.

Let $j \in \NN$. By \eqref{eq:thetanz}, \eqref{eq:Cjps} and the above density
result, there exist $f \in \cF$ and $\chi \in L^2(\Gamma_1)$ such that
$c_j' \neq 0$. Applying Lemma~\ref{lem:residue}(i) with this choice of $f$ and
$\chi$, we obtain that $-\lambda_j' \in \widetilde{\cP}$. This proves that
$\cP \subset \widetilde{\cP}$, and we have $\widetilde{\cP} \subset \cP$ as
well, upon exchanging the roles of $q$ and $\tq$. Therefore, the two strictly
increasing sequences $(\lambda_j')_{j \in \NN}$ and $(\tla_j')_{j \in \NN}$
enumerate the same set, so they coincide:
\begin{equation}\label{eq:eigeq}
\lambda_j' = \tla_j', \quad j \in \NN.
\end{equation}
Then, Lemma~\ref{lem:residue}(ii) and \eqref{eq:Cjps} yield
$\ps{\theta_j - \ttheta_j}{\overline{f} \otimes \chi}_{L^2(\Gamma_0 \times
\Gamma_1)} = 0$ for all $f \in \cF$ and all $\chi \in L^2(\Gamma_1)$, and hence,
by density,
\begin{equation}\label{eq:step4}
\theta_j = \ttheta_j \ \text{ a.e.\ in } \Gamma_0 \times \Gamma_1,
\quad j \in \NN.
\end{equation}

\subsection{Step 5: identification of the boundary spectral data}

The last step of the proof relies on the following algebraic result, which is
inspired by \cite{CK01} (see also \cite[Lemma~3.2]{Ch24b}).

\begin{lemma}\label{lem:alg}
Let $d,\, \widetilde{d} \in \NN$, and let $a = (a_1,\ldots,a_d)^T$ and
$\ta = (\ta_1,\ldots,\ta_{\widetilde{d}})^T$ be two vector-valued functions whose
components are real-valued functions in $L^2(\Gamma_0 \cup \Gamma_1)$. Assume
that for each $\Gamma' \in \{ \Gamma_0,\ \Gamma_1,\ \Gamma_0 \cap \Gamma_1 \}$,
the functions $a_1,\ldots,a_d$ (resp., $\ta_1,\ldots,\ta_{\widetilde{d}}$) are
linearly independent in $L^2(\Gamma')$, and that
\begin{equation}\label{eq:star}
a(\sigma) \cdot a(\sigma') = \ta(\sigma) \cdot \ta(\sigma')
\ \text{ for a.e. } (\sigma,\sigma') \in \Gamma_0 \times \Gamma_1.
\end{equation}
Then we have $d = \widetilde{d}$, and there exists an orthogonal matrix $O$ of
size $d$ with real entries such that
\[
a = O\, \ta \ \text{ a.e.\ in } \Gamma_0 \cup \Gamma_1.
\]
\end{lemma}

\begin{proof}
For $f \in L^2(\Gamma_0,\RR)$ we set
\[
a[f] := \int_{\Gamma_0} f(\sigma) a(\sigma)\, \dd\sigma \in \RR^d
\quad \text{and} \quad
\ta[f] := \int_{\Gamma_0} f(\sigma) \ta(\sigma)\, \dd\sigma \in
\RR^{\widetilde{d}}.
\]
The linear map $f \mapsto a[f]$ is onto $\RR^d$. Otherwise, there would exist
$c \in \RR^d \setminus \{0\}$ such that
$c \cdot a[f] = \int_{\Gamma_0} f\, (c \cdot a)\, \dd\sigma = 0$ for all
$f \in L^2(\Gamma_0,\RR)$, that is to say such that $c \cdot a = 0$ a.e.\ in
$\Gamma_0$, in contradiction with the linear independence of $a_1,\ldots,a_d$
in $L^2(\Gamma_0)$. Thus, we may find $f_1,\ldots,f_d$ in $L^2(\Gamma_0,\RR)$
such that $(a[f_1],\ldots,a[f_d])$ is the canonical basis of $\RR^d$.
Multiplying \eqref{eq:star} by $f_i(\sigma)$, $i = 1,\ldots,d$, and integrating
with respect to $\sigma$ over $\Gamma_0$, we get with the aid of Fubini's
theorem that
\[
a_i(\sigma') = a[f_i] \cdot a(\sigma') = \ta[f_i] \cdot \ta(\sigma')
\ \text{ for a.e. } \sigma' \in \Gamma_1,\ i = 1,\ldots,d.
\]
This can be equivalently rewritten as
\begin{equation}\label{eq:M1}
a = O_1\, \ta \ \text{ a.e.\ in } \Gamma_1,
\end{equation}
where $O_1$ is the $d \times \widetilde{d}$ matrix with real entries, whose
$i$-th row is $\ta[f_i]^T$, $i = 1,\ldots,d$. Since $a_1,\ldots,a_d$ are
linearly independent in $L^2(\Gamma_1)$, the rank of $O_1$ is equal to $d$, and
consequently we have $d \leq \widetilde{d}$. As $a$ and $\ta$ play symmetric
roles, we have $\widetilde{d} \leq d$ as well. Therefore, we have
$d = \widetilde{d}$ and $O_1$ is invertible.

Next, since \eqref{eq:star} is symmetric with respect to $\Gamma_0$ and
$\Gamma_1$, we obtain by arguing as above that
\begin{equation}\label{eq:M0}
a = O_0\, \ta \ \text{ a.e.\ in } \Gamma_0,
\end{equation}
for some invertible matrix $O_0$ of size $d$ with real entries. Now, inserting
\eqref{eq:M1} and \eqref{eq:M0} into \eqref{eq:star}, we find that
\[
\big( (O_1^T O_0 - I)\, \ta(\sigma) \big) \cdot \ta(\sigma') = 0
\ \text{ for a.e. } (\sigma,\sigma') \in \Gamma_0 \times \Gamma_1,
\]
where $I$ denotes the identity matrix of size $d$. Since
$\ta_1,\ldots,\ta_d$ are linearly independent in $L^2(\Gamma_1)$, this entails
that $(O_1^T O_0 - I)\, \ta = 0$ a.e.\ in $\Gamma_0$, and since they are
linearly independent in $L^2(\Gamma_0)$, we get that
\begin{equation}\label{eq:M1M0}
O_1^T O_0 = I.
\end{equation}
Finally, we deduce from \eqref{eq:M1} and \eqref{eq:M0} that
$(O_0 - O_1)\, \ta = 0$ a.e.\ in $\Gamma_0 \cap \Gamma_1$, and hence that
$O_0 = O_1$, since $\ta_1,\ldots,\ta_d$ are linearly independent in
$L^2(\Gamma_0 \cap \Gamma_1)$. Thus, the matrix $O := O_0 = O_1$ is orthogonal
by \eqref{eq:M1M0}, and the desired result follows from \eqref{eq:M1} and
\eqref{eq:M0}.
\end{proof}

\begin{proof}[End of the proof of Theorem~\ref{thm:bsd}]
Fix $j \in \NN$. Since $\Gamma_0$, $\Gamma_1$ and $\Gamma_0 \cap \Gamma_1$ are
nonempty open subsets of $\Gamma$, the components of $\Psi_j$ (resp., $\tPsi_j$)
are linearly independent in $L^2(\Gamma')$ for each
$\Gamma' \in \{ \Gamma_0,\ \Gamma_1,\ \Gamma_0 \cap \Gamma_1 \}$, by
Corollary~\ref{cor:indep}. In light of \eqref{eq:theta} and \eqref{eq:step4},
we may thus apply Lemma~\ref{lem:alg} with $a = \Psi_j$ and $\ta = \tPsi_j$.
We obtain that $m_j = \tm_j$ and that there exists an orthogonal matrix $O_j$ of
size $m_j$ with real entries, such that
\begin{equation}\label{eq:PsiM}
\Psi_j = O_j \tPsi_j \ \text{ a.e.\ in } \Gamma_0 \cup \Gamma_1.
\end{equation}
Since the distinct eigenvalues of $A$ and $\tA$ coincide by \eqref{eq:eigeq},
and have the same multiplicities, we have $\lambda_k = \tla_k$ for all
$k \in \NN$.

It remains to modify the eigenfunctions of $\tA$. We set
\[
(\tphi_{j,1}',\ldots,\tphi_{j,m_j}')^T
:= O_j\, (\tphi_{j,1},\ldots,\tphi_{j,m_j})^T.
\]
Since $O_j$ is an orthogonal matrix with real entries,
$\{ \tphi_{j,l}',\ l = 1,\ldots,m_j \}$ is a real-valued orthonormal basis of
$\ker(\tA - \tla_j')$, and by linearity of the trace operator and
\eqref{eq:PsiM}, we have
\[
\big( {\tphi_{j,1}'}{}_{|\Gamma},\ldots,{\tphi_{j,m_j}'}{}_{|\Gamma} \big)^T
= O_j \tPsi_j = \Psi_j \ \text{ a.e.\ in } \Gamma_0 \cup \Gamma_1.
\]
Therefore, the family $\{ \tphi_{j,l}',\ j \in \NN,\ l = 1,\ldots,m_j \}$,
relabelled as a sequence in accordance with $\{ \phi_k,\ k \in \NN \}$, is an
orthonormal basis of $H$ made of real-valued eigenfunctions of $\tA$, which
fulfills the claim of Theorem~\ref{thm:bsd}.
\end{proof}

\subsection{Proof of Theorem \ref{thm:uniq}}\label{sec:proof-uniq}

We start with the elliptic Robin problem associated with $A$ at a fixed energy.

\begin{lemma}\label{lem:ell}
Let $s \in [\lambda^*,\infty)$ and $f \in L^2(\Gamma)$. Then there exists a
unique $w \in V$ such that
\begin{equation}\label{eq:ell}
\fra(w,v) + s (w,v) = \ps{f}{v_{|\Gamma}}_\Gamma, \quad v \in V.
\end{equation}
Moreover, we have $s + \lambda_k \geq \mu_k > 0$ for all $k \in \NN$, and
\begin{equation}\label{eq:ellrep}
w = \sum_{k \geq 1} \frac{\ps{f}{\psi_k}_\Gamma}{s + \lambda_k}\, \phi_k
\ \text{ in } V,
\qquad
w_{|\Gamma} = \sum_{k \geq 1} \frac{\ps{f}{\psi_k}_\Gamma}{s + \lambda_k}\,
\psi_k \ \text{ in } L^2(\Gamma).
\end{equation}
\end{lemma}

\begin{proof}
Since $s \geq \lambda^*$, we infer from \eqref{eq:acont} and \eqref{eq:coer}
that
\begin{equation}\label{eq:ellcoer}
\kappa \norm{u}_V^2 \leq \fra(u,u) + s \norm{u}_H^2
\leq (C_\fra + s) \norm{u}_V^2, \quad u \in V,
\end{equation}
so the Hermitian form $\fra + s(\cdot,\cdot)$ is a scalar product on $V$, whose
norm is equivalent to the usual one. As the antilinear form
$v \mapsto \ps{f}{v_{|\Gamma}}_\Gamma$ is continuous on $V$ by \eqref{eq:trace},
the existence and the uniqueness of $w$ follow from the Riesz representation
theorem. Next, we have $s + \lambda_k \geq \lambda^* + \lambda_k = \mu_k > 0$.
Further, taking $v = \phi_k$ in \eqref{eq:ell}, we get from \eqref{eq:eig}, upon
recalling that $\fra$ is Hermitian and that $\lambda_k$ and $\psi_k$ are
real-valued, that
\[
(s + \lambda_k)(w,\phi_k) = \overline{\fra(\phi_k,w) + s(\phi_k,w)}
= \fra(w,\phi_k) + s(w,\phi_k) = \ps{f}{\psi_k}_\Gamma, \quad k \in \NN.
\]
Since $w \in V$, the first identity of \eqref{eq:ellrep} follows from this and
Lemma~\ref{lem:spec}(i), and the second one is then a consequence of the
continuity of the trace operator from $V$ into $L^2(\Gamma)$.
\end{proof}

In the remaining part of this subsection, we assume that
$\Gamma_0 \cup \Gamma_1 = \Gamma$ and that \eqref{eq:hyp} holds; this is the case
under the assumptions of Theorem~\ref{thm:uniq}. Since the proof of
Theorem~\ref{thm:bsd} given above only relies on \eqref{eq:hyp}, we may choose
the orthonormal basis $\{ \tphi_k,\ k \in \NN \}$ of eigenfunctions of $\tA$ in
such a way that
\begin{equation}\label{eq:bsdfull}
\lambda_k = \tla_k \ \text{ and } \ \psi_k = \tpsi_k \ \text{ a.e.\ on } \Gamma,
\quad k \in \NN.
\end{equation}
We denote by $\widetilde{\fra}$ the sesquilinear form defined by \eqref{eq:form}
with $\tq$ in place of $q$, and we set
\[
\fra_\sharp(u,v) := \fra(u,v) + \lambda_\sharp (u,v)
\quad \text{and} \quad
\widetilde{\fra}_\sharp(u,v) := \widetilde{\fra}(u,v) + \lambda_\sharp (u,v),
\quad u,\, v \in V.
\]
As $\lambda_\sharp \geq \lambda^*$ and $\lambda_\sharp \geq \tla^*$, we may apply
Lemma~\ref{lem:ell} with $s = \lambda_\sharp$, both to $q$ and to $\tq$. For
$f \in L^2(\Gamma)$, we denote by $w_f$ (resp., $\widetilde{w}_f$) the unique
function $w \in V$ such that $\fra_\sharp(w,v) = \ps{f}{v_{|\Gamma}}_\Gamma$
(resp., $\widetilde{\fra}_\sharp(w,v) = \ps{f}{v_{|\Gamma}}_\Gamma$) for all
$v \in V$. We recall from \eqref{eq:ellcoer} that $\fra_\sharp$ and
$\widetilde{\fra}_\sharp$ are scalar products on $V$, whose norms are equivalent
to the usual one. Finally, we introduce the spaces of solutions
\begin{align*}
\mathcal{E} & := \{ u \in V;\ -\Delta u + (q + \lambda_\sharp) u = 0
\ \text{ in } \mathcal{D}'(\Omega) \}, \\
\widetilde{\mathcal{E}} & := \{ u \in V;\ -\Delta u + (\tq + \lambda_\sharp) u = 0
\ \text{ in } \mathcal{D}'(\Omega) \}.
\end{align*}
Notice that $(q + \lambda_\sharp) u \in L^p(\Omega)$ for all $u \in V$, by
H\"older's inequality, since $V \subset L^{p^*}(\Omega)$ and
$1/p = 2/n + 1/p^*$. Moreover, as $q$, $\tq$ and $\lambda_\sharp$ are real, we
have $\overline{u} \in \mathcal{E}$ (resp., $\overline{u} \in
\widetilde{\mathcal{E}}$) whenever $u \in \mathcal{E}$ (resp., $u \in
\widetilde{\mathcal{E}}$).

\begin{lemma}\label{lem:dense}
\begin{enumerate}
\item[\textup{(i)}] A function $u \in V$ belongs to $\mathcal{E}$ if and only if
$\fra_\sharp(u,v) = 0$ for all $v \in H_0^1(\Omega)$.
\item[\textup{(ii)}] The set $\mathcal{W} := \{ w_f;\ f \in L^2(\Gamma) \}$ is a
dense subspace of $\mathcal{E}$, endowed with the norm of $V$.
\end{enumerate}
The same is true for $\widetilde{\mathcal{E}}$, $\widetilde{\fra}_\sharp$ and the
$\widetilde{w}_f$, $f \in L^2(\Gamma)$.
\end{lemma}

\begin{proof}
(i) Let $u \in V$ and $\varphi \in C_c^\infty(\Omega)$. As the trace of
$\overline{\varphi}$ on $\Gamma$ vanishes, we have
\[
\fra_\sharp(u,\overline{\varphi}) = \int_\Omega \nabla u \cdot \nabla \varphi\,
\dd x + \int_\Omega (q + \lambda_\sharp)\, u\, \varphi\, \dd x
= \big\langle -\Delta u + (q + \lambda_\sharp) u , \varphi \big\rangle,
\]
where $\langle \cdot,\cdot \rangle$ denotes the duality pairing between
$\mathcal{D}'(\Omega)$ and $C_c^\infty(\Omega)$. Thus $u \in \mathcal{E}$ if and
only if $\fra_\sharp(u,v) = 0$ for all $v \in C_c^\infty(\Omega)$. Since
$C_c^\infty(\Omega)$ is dense in $H_0^1(\Omega)$ and the antilinear form
$\fra_\sharp(u,\cdot)$ is continuous on $V$, this proves (i).

(ii) By (i), $\mathcal{E}$ is the intersection of the kernels of the continuous
linear forms $u \mapsto \fra_\sharp(u,v)$, $v \in H_0^1(\Omega)$, so it is a
closed subspace of $V$. Next, $\mathcal{W}$ is a linear subspace of $V$, since
$w_f$ depends linearly on $f$, and we have $\mathcal{W} \subset \mathcal{E}$ by
(i), since $\fra_\sharp(w_f,v) = \ps{f}{v_{|\Gamma}}_\Gamma = 0$ for all
$v \in H_0^1(\Omega)$. Now, let $u \in \mathcal{E}$ satisfy
$\fra_\sharp(u,w_f) = 0$ for all $f \in L^2(\Gamma)$. Then we have
\[
\ps{f}{u_{|\Gamma}}_\Gamma = \fra_\sharp(w_f,u) = \overline{\fra_\sharp(u,w_f)}
= 0, \quad f \in L^2(\Gamma),
\]
and taking $f = u_{|\Gamma}$ yields $u_{|\Gamma} = 0$. Since $\Gamma$ is of
class $C^1$, this entails that $u \in H_0^1(\Omega)$, see
\cite[Section~5.5, Theorem~2]{Eva10}, and hence $\fra_\sharp(u,u) = 0$ by (i).
Therefore, $u = 0$. This shows that the orthogonal complement of $\mathcal{W}$
in the Hilbert space $\mathcal{E}$, endowed with the scalar product
$\fra_\sharp$, is reduced to $\{0\}$, so that $\mathcal{W}$ is dense in
$\mathcal{E}$ for the norm associated with $\fra_\sharp$, and hence for the norm
of $V$. The case of $\widetilde{\mathcal{E}}$ is treated in the same way.
\end{proof}

\begin{lemma}\label{lem:identity}
Assume \eqref{eq:bsdfull}. Then we have
\begin{equation}\label{eq:identity}
\int_\Omega (q - \tq)\, u\, v\, \dd x = 0, \quad u \in \mathcal{E},\
v \in \widetilde{\mathcal{E}}.
\end{equation}
\end{lemma}

\begin{proof}
For all $u,\, v \in V$, we have $(q - \tq) u v \in L^1(\Omega)$ and
\begin{equation}\label{eq:holder3}
\Big\lvert \int_\Omega (q - \tq)\, u\, v\, \dd x \Big\rvert
\leq \norm{q - \tq}_{L^{n/2}(\Omega)} \norm{u}_{L^{p^*}(\Omega)}
\norm{v}_{L^{p^*}(\Omega)},
\end{equation}
by H\"older's inequality, since $2/n + 2/p^* = 1$. Let $f \in L^2(\Gamma)$. By
the second identity of \eqref{eq:ellrep}, applied with $s = \lambda_\sharp$ to $q$
and to $\tq$, and by \eqref{eq:bsdfull}, we have
${w_f}_{|\Gamma} = {\widetilde{w}_f}{}_{|\Gamma}$, and consequently
$w_f - \widetilde{w}_f \in H_0^1(\Omega)$ by
\cite[Section~5.5, Theorem~2]{Eva10}. Let $v \in \widetilde{\mathcal{E}}$, so
that $\overline{v} \in \widetilde{\mathcal{E}}$. By definition of $w_f$ and
$\widetilde{w}_f$, we have
\[
\fra_\sharp(w_f,\overline{v}) = \ps{f}{\overline{v}_{|\Gamma}}_\Gamma
= \widetilde{\fra}_\sharp(\widetilde{w}_f,\overline{v}),
\]
while \eqref{eq:form} yields
$\fra_\sharp(w_f,\overline{v}) - \widetilde{\fra}_\sharp(w_f,\overline{v})
= \int_\Omega (q - \tq)\, w_f\, v\, \dd x$. Therefore, we get
\[
\int_\Omega (q - \tq)\, w_f\, v\, \dd x
= \widetilde{\fra}_\sharp(\widetilde{w}_f - w_f,\overline{v})
= \overline{\widetilde{\fra}_\sharp(\overline{v},\widetilde{w}_f - w_f)} = 0,
\]
by Lemma~\ref{lem:dense}(i), applied to $\tq$. Finally, let $u \in \mathcal{E}$.
By Lemma~\ref{lem:dense}(ii), there exists a sequence $(f_j)_{j \in \NN}$ in
$L^2(\Gamma)$ such that $w_{f_j} \to u$ in $V$ as $j \to \infty$, hence in
$L^{p^*}(\Omega)$. We deduce from this and \eqref{eq:holder3} that
$\int_\Omega (q - \tq) u v\, \dd x = \lim_{j \to \infty} \int_\Omega (q - \tq)
w_{f_j} v\, \dd x = 0$.
\end{proof}

\begin{proof}[End of the proof of Theorem~\ref{thm:uniq}]
We use complex geometric optics solutions. Let $B$ be an open ball of $\RR^n$
such that $\overline{\Omega} \subset B$, and put
$q_1 := (q + \lambda_\sharp) \mathbf{1}_\Omega$ and
$q_2 := (\tq + \lambda_\sharp) \mathbf{1}_\Omega$, where $\mathbf{1}_\Omega$ is
the characteristic function of $\Omega$. Then $q_1$ and $q_2$ belong to
$L^{n/2}(B)$, and they vanish outside $\Omega$. Since the boundary of $B$ is
smooth, we may apply \cite[Proposition~4.1]{KU16} with $m = 1$, upon
substituting $B$ for $\Omega$ and $q_j$ for $q$, $j = 1,2$. We obtain that for
every $\zeta \in \CC^n$ such that $\zeta \cdot \zeta = 0$ and $\abs{\zeta}$ is
sufficiently large, there exists a solution to the equation
$-\Delta u + q_j u = 0$ in the distributional sense in $\RR^n$, of the form
$u(x) = e^{x \cdot \zeta} (1 + r(x))$, where the remainder $r$ satisfies
\begin{equation}\label{eq:cgo}
\norm{r}_{L^{p^*}(\RR^n)} = O(1)
\quad \text{and} \quad
\norm{r}_{L^2(K)} \to 0
\quad \text{as } \abs{\zeta} \to \infty,
\end{equation}
for every compact subset $K$ of $\RR^n$. Moreover, $u_{|B} \in H^1(B)$ by
\cite[Remark~4.2]{KU16}, so that $u_{|\Omega} \in V$. As $q_1 = q + \lambda_\sharp$
and $q_2 = \tq + \lambda_\sharp$ in $\Omega$, we thus have
$u_{|\Omega} \in \mathcal{E}$ when $j = 1$, and
$u_{|\Omega} \in \widetilde{\mathcal{E}}$ when $j = 2$.

Fix $\xi \in \RR^n$. As $n \geq 3$, there exist two unit vectors $\eta_1$ and
$\eta_2$ of $\RR^n$ such that $\eta_1 \cdot \eta_2 = \eta_1 \cdot \xi =
\eta_2 \cdot \xi = 0$. For $\tau \in [\abs{\xi}/2,\infty)$, we set
\[
\zeta_1 := \tau \eta_1 + i \Big( \frac{\xi}{2}
+ \Big( \tau^2 - \frac{\abs{\xi}^2}{4} \Big)^{1/2} \eta_2 \Big),
\qquad
\zeta_2 := -\tau \eta_1 + i \Big( \frac{\xi}{2}
- \Big( \tau^2 - \frac{\abs{\xi}^2}{4} \Big)^{1/2} \eta_2 \Big),
\]
in such a way that $\zeta_j \cdot \zeta_j = 0$ and $\abs{\zeta_j} = \sqrt{2}\,
\tau$ for $j = 1,2$, and $\zeta_1 + \zeta_2 = i \xi$. Hence, for $\tau$
sufficiently large, there exist $u_1 = e^{x \cdot \zeta_1} (1 + r_1) \in
\mathcal{E}$ and $u_2 = e^{x \cdot \zeta_2} (1 + r_2) \in \widetilde{\mathcal{E}}$
(we still denote by $u_j$ the restriction to $\Omega$ of the above solutions),
where $r_1$ and $r_2$ satisfy \eqref{eq:cgo} as $\tau \to \infty$. Putting
$Q := q - \tq \in L^{n/2}(\Omega,\RR)$, we infer from \eqref{eq:identity} that
\begin{equation}\label{eq:fourier}
\int_\Omega Q(x)\, e^{i \xi \cdot x}\, \dd x
= - \int_\Omega Q(x)\, e^{i \xi \cdot x} \big( r_1 + r_2 + r_1 r_2 \big)(x)\,
\dd x.
\end{equation}
Let $\varepsilon > 0$. By the dominated convergence theorem, there exists
$N > 0$ such that $Q^\sharp := Q\, \mathbf{1}_{\{\abs{Q} \leq N\}}$ satisfies
$\norm{Q - Q^\sharp}_{L^{n/2}(\Omega)} \leq \varepsilon$. By the
Cauchy--Schwarz inequality, we have
\[
\Big\lvert \int_\Omega Q^\sharp e^{i \xi \cdot x} \big( r_1 + r_2 + r_1 r_2
\big)\, \dd x \Big\rvert
\leq N \Big( \abs{\Omega}^{1/2} \big( \norm{r_1}_{L^2(\Omega)}
+ \norm{r_2}_{L^2(\Omega)} \big) + \norm{r_1}_{L^2(\Omega)}
\norm{r_2}_{L^2(\Omega)} \Big),
\]
and the right-hand side tends to zero as $\tau \to \infty$, by \eqref{eq:cgo}
with $K = \overline{\Omega}$. Similarly, since $2/n + 2/p^* = 1$, H\"older's
inequality yields
\begin{align*}
\Big\lvert \int_\Omega (Q - Q^\sharp) e^{i \xi \cdot x} \big( r_1 + r_2
+ r_1 r_2 \big)\, \dd x \Big\rvert
& \leq \varepsilon \Big( \abs{\Omega}^{1/p^*} \big( \norm{r_1}_{L^{p^*}(\RR^n)}
+ \norm{r_2}_{L^{p^*}(\RR^n)} \big) \\
& \qquad + \norm{r_1}_{L^{p^*}(\RR^n)} \norm{r_2}_{L^{p^*}(\RR^n)} \Big),
\end{align*}
and the right-hand side is bounded by $C \varepsilon$ for all $\tau$
sufficiently large, where $C$ is a positive constant independent of $\tau$ and
$\varepsilon$, by \eqref{eq:cgo}. Thus, sending $\tau$ to infinity in
\eqref{eq:fourier}, we get that
$\abs{\int_\Omega Q(x) e^{i \xi \cdot x}\, \dd x} \leq C \varepsilon$, and hence
$\int_\Omega Q(x) e^{i \xi \cdot x}\, \dd x = 0$, as $\varepsilon$ is
arbitrary. Since $\xi \in \RR^n$ is arbitrary, the Fourier transform of the
extension of $Q$ by zero outside $\Omega$, which belongs to $L^1(\RR^n)$,
vanishes identically. Therefore, we have $Q = 0$, i.e.\ $q = \tq$ in $\Omega$. This
completes the proof of Theorem~\ref{thm:uniq}. Finally, since the above proofs
of Theorems~\ref{thm:uniq} and \ref{thm:bsd} only rely on the assumption
\eqref{eq:hyp}, the claim of Remark~\ref{rem:main}(ii) is established as well.
\end{proof}

\begin{remark}\label{rem:elliptic}
The meromorphic function $M$ of Lemma~\ref{lem:mero} has a simple
interpretation in terms of elliptic boundary value problems. Let
$s \in [\lambda^*,\infty)$, let $f \in L^2(\Gamma)$ with $\supp f \subset
\Gamma_0$, and let $w_s \in V$ be the solution to \eqref{eq:ell} given by
Lemma~\ref{lem:ell}, i.e.\ the weak solution to $(-\Delta + q + s) w = 0$ in
$\Omega$ and $(\partial_\nu + \alpha) w = f$ on $\Gamma$. Then it follows from
\eqref{eq:ellrep} and \eqref{eq:ck} that
\[
\ps{{w_s}_{|\Gamma_1}}{\chi}_{\Gamma_1} = M(s), \quad \chi \in L^2(\Gamma_1).
\]
Moreover, the map $f \mapsto {w_s}_{|\Gamma_1}$ is bounded from $L^2(\Gamma)$
into $L^2(\Gamma_1)$, since $\kappa \norm{w_s}_V^2 \leq \frn
\norm{f}_{L^2(\Gamma)} \norm{w_s}_V$ by \eqref{eq:ellcoer}. Therefore,
Steps~1 to 3 of the proof of Theorem~\ref{thm:bsd}, together with the density
in $L^2(\Gamma_0)$ of the set of functions $f \in L^2(\Gamma)$ such that
$\supp f \subset \Gamma_0$, show that the Robin-to-Dirichlet map $D_q$
uniquely determines the partial elliptic Robin-to-Dirichlet maps
$f \mapsto {w_s}_{|\Gamma_1}$, with $f \in L^2(\Gamma_0)$, for every
$s \in (\lambda_\sharp,\infty)$. When $\Gamma_0 \cup \Gamma_1 = \Gamma$,
Theorem~\ref{thm:bsd} and Lemma~\ref{lem:ell} show that $D_q$ even determines the
full elliptic Robin-to-Dirichlet map $f \in L^2(\Gamma) \mapsto
{w_s}_{|\Gamma} \in L^2(\Gamma)$, which is the starting point of the above
proof of Theorem~\ref{thm:uniq}.
\end{remark}


\end{document}